\documentclass[11pt,twoside]{amsart} 
\usepackage{amssymb}
\advance \topmargin by -\headheight
\advance \topmargin by -\headsep
\evensidemargin \oddsidemargin
\newtheorem{theorem}{Theorem}[section]

\newtheorem{corollary}[theorem]{Corollary}
\newtheorem{lemma}[theorem]{Lemma}
\newtheorem{proposition}[theorem]{Proposition}

\title{Admissible semi-linear representations} 
\author{M.Rovinsky} 
\thanks{The author was supported in part by RFBR grant 02-01-22005}
\begin{document} 
\begin{abstract} The category of admissible (in the appropriately
modified sense of representation theory of totally disconnected 
groups) semi-linear representations of the automorphism group of 
an algebraically closed extension of infinite transcendence degree 
of the field of algebraic complex numbers is described. \end{abstract}
\maketitle 

Let $k$ be a field of characteristic zero containing all $\ell$-primary 
roots of unity for a prime $\ell$, $F$ be a universal domain over $k$, 
i.e., an algebraically closed extension of $k$ of countable 
transcendence degree, and $G_{F/k}$ be the field automorphism group 
of $F$ over $k$. We consider $G_{F/k}$ as a topological group with 
the base of open subgroups generated by $\{G_{F/k(x)}~|~x\in F\}$. 

Denote by ${\mathcal C}$ the category of smooth (with open stabilizers) 
$F$-semi-linear representations of $G_{F/k}$, i.e., $F$-vector spaces 
$V$ endowed with an additive semi-linear ($g(fv)=gf\cdot gv$ for any 
$f\in F$, $g\in G_{F/k}$ and $v\in V$) action $G_{F/k}\times V\to V$ 
of $G_{F/k}$. 

Denote by ${\mathcal A}$ the full sub-category of ${\mathcal C}$ 
whose objects $V$ are {\sl admissible}: $\dim_{F^U}V^U<\infty$ for 
any open subgroup $U\subseteq G_{F/k}$. Clearly, ${\mathcal A}$ is 
an additive category and it is shown in \cite{pgl} that it is a 
tensor (but not rigid) category. In the present paper one proves 
that the category ${\mathcal A}$ is abelian (Theorem \ref{final}), 
and $F$ is its projective object (Proposition \ref{proj-F}).  

Let the ideal ${\mathfrak m}\subset F\otimes_{k_0}F$ be the kernel of the 
multiplication map $F\otimes_{k_0}F\stackrel{\times}{\longrightarrow}F$, 
where $k_0=k\cap\overline{{\mathbb Q}}$ is the number subfield of $k$. 
Consider the powers ${\mathfrak m}^s\subseteq F\otimes_{k_0}F$ 
of the ideal ${\mathfrak m}$ for all $s\ge 0$ as objects of 
${\mathcal C}$ with the $F$-multiplication via $F\otimes_{k_0}k_0$. 

In this paper we study the category ${\mathcal A}$ and 
describe it if $k$ is a number field. Namely, in the 
case $k=\overline{{\mathbb Q}}$ we prove the following: 
\begin{itemize} \item The sum of the images of the $F$-tensor powers 
$\bigotimes^{\ge\bullet}_F{\mathfrak m}$ under all morphisms in 
${\mathcal C}$ defines a decreasing filtration $W^{\bullet}$ on the 
objects of ${\mathcal A}$ such that its graded quotients $gr^q_W$ are 
finite direct sums of direct summands of $\bigotimes^q_F\Omega^1_F$ 
(cf. \S\ref{ves-filtr}, p.\pageref{ves-filtr} and Theorem 
\ref{irr-form}). This filtration is evidently functorial 
and multiplicative: $(W^pV_1)\otimes_F(W^qV_2)\subseteq 
W^{p+q}(V_1\otimes_FV_2)$ for any $p,q\ge 0$ and any 
$V_1,V_2\in{\mathcal A}$. 
\item ${\mathcal A}$ is equivalent to the direct sum of the category 
of finite-dimensional $k$-vector spaces and its abelian full 
subcategory ${\mathcal A}^{\circ}$ with objects 
$V$ such that $V^{G_{F/k}}=0$ (Lemma \ref{A-spl}). 
\item Any object $V$ of ${\mathcal A}^{\circ}$ is a quotient of a direct 
sum of objects (of finite length) of type $\bigotimes^q_F({\mathfrak m}
/{\mathfrak m}^s)$ for some $q,s\ge 1$ (Theorem \ref{irr-form}). 
\item If $V\in{\mathcal A}$ is of finite type then it is 
of finite length and $\dim_k{\rm Ext}_{{\mathcal A}}^j(V,V')
<\infty$ for any $j\ge 0$ and any $V'\in{\mathcal A}$; 
if $V\in{\mathcal A}$ is irreducible and ${\rm Ext}^1_{{\mathcal A}}
({\mathfrak m}/{\mathfrak m}^q,V)\neq 0$ for some $q\ge 2$ then 
$V\cong{\rm Sym}^q_F\Omega^1_F$ and ${\rm Ext}^1_{{\mathcal A}}
({\mathfrak m}/{\mathfrak m}^q,V)\cong k$ (Corollary \ref{fin-dim-ext}). 
\item ${\mathcal A}^{\circ}$ has no projective objects
(Corollary \ref{edinstv-proekt}), but $\bigotimes^q_F{\mathfrak m}$ 
are its ``projective pro-generators'': the functor 
${\rm Hom}_{{\mathcal C}}(\bigotimes^q_F{\mathfrak m},-)$ is exact 
on ${\mathcal A}$ for any $q$ (Corollary \ref{4.12}). \end{itemize} 

\vspace{5mm}

To describe the objects of ${\mathcal A}$, one studies 
first their ``restrictions'' to projective groups 
($\cong{\rm PGL}_mk$), considered as subquotients of $G_{F/k}$. 
It is known (\cite{pgl}) that such semi-linear representations 
are related to homogeneous vector bundles on projective spaces. 

Let $n\ge 1$ be an integer, $K_n=k({\mathbb P}^n_k)$ be the function 
field of an $n$-dimensional projective $k$-space ${\mathbb P}^n_k$ and 
$G_n={\rm Aut}({\mathbb P}^n_k/k)$ be its automorphism group. 

Fix a $k$-field embedding $K_n\hookrightarrow F$. We show, 
in particular, that if $V$ is an admissible $F$-semi-linear 
representation of $G_{F/k}$ with no sub-objects isomorphic to 
$F$ then any irreducible subquotient of the $K_n$-semi-linear 
representation $V^{G_{F/K_n}}$ of $G_n$ is a direct summand of 
$\bigotimes_{K_n}^{\ge 1}\Omega^1_{K_n/k}$ (this is Theorem 
\ref{main} and Proposition \ref{proj-F} below).

\subsection{Some motivation} The study of semi-linear representations    
comes from the study of ${\mathbb Q}$-linear representations of 
$G_{F/k}$, that are related to geometry, cf. \cite{pgl}. 

Let ${\mathcal S}m_G$ be the category of smooth representations of 
$G_{F/k}$ over $k$. Extending of coefficients to $F$ gives a faithful 
functor $F\otimes_k:{\mathcal S}m_G\longrightarrow{\mathcal C}$. 
It is not full: if $U\subset G_{F/k}$ is an open subgroup and
$\overline{f}\in(F^{\times}/k^{\times})^U-\{1\}$ then $[\sigma]
\mapsto\sigma f\cdot[\sigma]$ defines an element of ${\rm End}
_{{\mathcal C}}(F[G_{F/k}/U])$ which is not in ${\rm End}
_{{\mathcal S}m_G}(k[G_{F/k}/U])$. However, its restriction 
to the subcategory ${\mathcal I}_G\otimes k$ of ``homotopy 
invariant'' representations\footnote{i.e. such that 
$W^{G_{F/L}}=W^{G_{F/L'}}$ for any purely transcendental 
extension $L'/L$ of subfields in $F$} is. 
\begin{lemma} If $k=\overline{k}$ then the functor 
${\mathcal I}_G\otimes k\stackrel{F\otimes_k}{\longrightarrow}
{\mathcal C}$ is fully faithful. \end{lemma}
{\it Proof.} More generally, let us show that ${\rm Hom}_{{\mathcal S}m_G}
(W,W')={\rm Hom}_{{\mathcal C}}(F\otimes_kW,F\otimes_kW')$ for any 
$W\in{\mathcal I}_G\otimes k$ and any $W'\in{\mathcal S}m_G$. Let 
$\varphi\in{\rm Hom}_G(W,F\otimes_kW')$ and $\varphi(w)=\sum_{j=1}^N
f_j\otimes w_j$ for some $w\in W$, $w_j\in W'$, $f_j\in F$ 
and minimal possible $N\ge 1$. We have to show that $f_j\in k$. 

Choose a smooth proper model $X$ of $k(f_1,\dots,f_N)$ over $k$. 
If it is not a point, choose a generically finite rational 
dominant map $\pi$ to a projective space $Y$ over $k$ which is 
\begin{itemize} \item well-defined at the generic points of the irreducible 
components $D_{\alpha}$ of the divisors of poles of $f_1,\dots,f_N$, 
\item induces on each $D_{\alpha}$ a birational map and 
\item separates $D_{\alpha}$. \end{itemize} Then the 
trace $\pi_{\ast}\varphi(w)$ has poles. On the other hand,
$\pi_{\ast}\varphi(w)$ is in the image of $W^{G_{F/k(Y)}}=
W^{G_{F/k}}$, so $\pi_{\ast}\varphi(w)\in(F\otimes_kW')^{G_{F/k}}
=k\otimes_k(W')^{G_{F/k}}$ by Lemma 7.5 of \cite{pgl}. 
This contradiction implies that $f_j\in k$, and therefore, 
$\varphi(W)\subseteq k\otimes_kW'\subseteq F\otimes_kW'$. \qed 

\vspace{4mm}

The ${\mathbb Q}$-linear representations of $G_{F/k}$ of particular 
interest are admissible representations, forming a full subcategory 
in ${\mathcal I}_G\otimes k$. Though tensoring with $F$ does not 
transform them to admissible semi-linear representations,\footnote{and 
moreover, there are irreducible objects of ${\mathcal C}$ outside of 
${\mathcal A}$ (Corollary \ref{exam-non-adm}).} there exists, 
at least if $k=\overline{{\mathbb Q}}$, a similar faithful 
functor in the opposite direction.

Namely, it is explained in Corollary \ref{funk-A-presh} that, for 
any object $V$ of ${\mathcal A}$ and any smooth $k$-variety $Y$, 
embedding of the generic points of $Y$ into $F$ determines a locally 
free coherent sheaf ${\mathcal V}_Y$ on $Y$. Any dominant morphism 
$X\stackrel{\pi}{\longrightarrow}Y$ of smooth $k$-varieties induces 
an injection of coherent sheaves $\pi^{\ast}{\mathcal V}_Y\hookrightarrow
{\mathcal V}_X$, which is an isomorphism if $\pi$ is {\'e}tale.

This gives an equivalence ${\mathcal S}:{\mathcal A}\stackrel{\sim}
{\longrightarrow}\{\mbox{``coherent'' sheaves in the smooth topology}\}$, 
$V\longmapsto(Y\mapsto{\mathcal V}_Y(Y))$. More generally, the
``coherent'' sheaves are contained in the category ${\mathcal F}l$ 
of the flat ``quasi-coherent'' sheaves in the smooth topology, 
cf. \S\ref{smooth-topol}, p.\pageref{smooth-topol}. For any 
flat ``quasi-coherent'' sheaf ${\mathcal V}$ in the smooth topology 
the space $\Gamma(Y,{\mathcal V}_Y)$ is a birational invariant 
of a proper $Y$ (Lemma \ref{birac-invar}). Then we get a left exact 
functor ${\mathcal F}l\stackrel{\Gamma}{\longrightarrow}
{\mathcal S}m_G$ given by $V\mapsto\lim\limits_{\longrightarrow}
\Gamma(Y,{\mathcal V}_Y)$, where $Y$ runs over the smooth proper 
models of subfields in $F$ of finite type over $k$. 

The functor $\Gamma\circ{\mathcal S}$ is faithful, since 
$\Gamma(Y',{\mathcal V}_{Y'})$ generates the (generic fibre of 
the) sheaf ${\mathcal V}_{Y'}$ for appropriate finite covers $Y'$ of $Y$ 
(Lemma \ref{birac-invar}), if ${\mathcal V}$ is ``coherent''. But it is 
not full, and the objects in its image are highly reducible. If $\Gamma
(Y,{\mathcal V}_Y)$ has the Galois descent property then $\Gamma(V)$ is 
admissible. However, there is no Galois descent property in general. 

\subsection{Notation} \label{obozn} Let $k$, $F$ and $G_{F/k}$ be 
as above. For a subfield $L$ of $F$ we denote by $\overline{L}$ 
its algebraic closure in $F$. We fix a transcendence basis 
$x_1,x_2,x_3,\dots$ of $F$ over $k$. 

For each $n\ge 1$ set $Y_n={\bf Spec}k[x_1^{\pm 1},\dots,x_n^{\pm 1}]
\subset{\mathbb A}^n_k={\bf Spec}k[x_1,\dots,x_n]\subset
{\mathbb P}^n_k={\bf Proj}k[X_0,\dots,X_n]$ with $x_j=X_j/X_0$, 
$K_n=k(x_1,\dots,x_n)$, $G_n={\rm Aut}({\mathbb P}^n_k/k)$. 

Let ${\rm Aff}_n=G_n\cap{\rm Aut}({\mathbb A}^n_k/k)$ be the 
affine subgroup, and $({\rm Aff}_n)_u$ be its unipotent radical, 
i.e., the translation subgroup. Let $H_n={\rm Aff}_n\cap
{\rm Aut}({\mathbb A}^n_k/{\mathbb A}^{n-1}_k)$ be the subgroup 
fixing the coordinates $x_1,\dots,x_{n-1}$ on ${\mathbb A}^n_k$. 
Let $T_n\subset G_n$ be the maximal torus acting freely on $Y_n$. 

Denote by $T_n^{{\rm tors}}$ the torsion subgroup in $T_n$. 

For a field extension $L/K$ we denote by ${\rm Der}(L/K)$ the Lie 
$K$-algebra of derivations of $L$ over $K$. For an integer $\ell\ge 2$, 
the group of $\ell$-th roots of unity in $\overline{k}$ is denoted by 
$\mu_{\ell}$, and the corresponding cyclotomic number subfield in 
$\overline{k}$ is denoted by ${\mathbb Q}(\mu_{\ell})$. 

\subsection{Structure of the paper} As it is mentioned above, we consider 
the projective groups $G_n$ as subquotients of $G_{F/k}$. In \S\ref{no-1} 
we identify irreducible subquotients of ``restrictions'' of objects of 
${\mathcal A}$ to $G_n$ with the generic fibres of the $G_n$-equivariant 
coherent sheaves on ${\mathbb P}^n_k$. Main ingredients there come from 
\cite{bt} and \cite{pgl}. In \S\ref{polozhit} we exclude some cases, 
thus showing that these irreducible subquotients are direct summands of 
$\bigotimes^{\bullet}_{K_n}\Omega^1_{K_n/k}$. In \S\ref{rasshir-v-A} we 
show that ${\mathcal A}$ is abelian and calculate ${\rm Ext}^1$-groups 
between the irreducible objects of a tannakian category 
$\mathfrak{SL}^n_u$ (defined at the beginning of \S\ref{no-1}, 
p.\pageref{no-1}) of semi-linear representations of $G_n$, containing 
``restrictions'' of objects of ${\mathcal A}$ to $G_n$. The latter 
part uses \cite{lr}. After showing principal structural results on 
${\mathcal A}$ (in \S\ref{str-ra-A-Q}) we identify (in \S\ref{smooth-topol}) 
${\mathcal A}$ with the category of ``coherent'' sheaves in smooth topology. 
Finally (in \S\ref{faktor-kat}), we define a descending filtration 
${\mathcal A}_{>\bullet}$ of ${\mathcal A}$ by Serre ``ideal'' subcategories. 
Then we localize the quotients ${\mathcal A}/{\mathcal A}_{>m}$ for each 
$m\ge 0$ to get a tannakian subcategory of finite-dimensional semi-linear 
representations of $G_{F'/k}$ over $F'$ for an algebraically closed 
extension $F'$ of $k$ in $F$ of transcendence degree $m$. 

\section{Equivariantness of irreducible $PGL$-sheaves} \label{no-1} 
Let $\mathfrak{SL}^u_n$ be the category of finite-dimensional 
semi-linear representations of $G_n$ over $K_n$ whose restrictions to 
the maximal torus $T_n$ in $G_n$ are of type $K_n\otimes_kW$ for unipotent 
representations $W$ of $T_n$ (where $T_n$ is considered as a discrete group). 

Note that $V=V^{T_n^{{\rm tors}}}\otimes_kK_n$ 
for any $V\in\mathfrak{SL}^u_n$. 

In \cite{pgl}, for $n>1$, a fully faithful functor $\mathfrak{SL}^u_n
\stackrel{{\mathcal S}}{\to}\{\mbox{coherent $G_n$-sheaves on 
${\mathbb P}^n_k$}\}$ is constructed. (A $G_n$-sheaf is $G_n$-equivariant 
sheaf if $G_n$ is considered as a discrete group. In other words, 
${\mathcal V}$ is a $G_n$-sheaf if it is endowed with a collection of 
isomorphisms $\alpha_g:{\mathcal V}\stackrel{\sim}{\longrightarrow} 
g^{\ast}{\mathcal V}$ for each $g\in G_n$ satisfying the chain rule: 
$\alpha_{hg}=g^{\ast}\alpha_h\circ\alpha_g$ for any $g,h\in G_n$. 
The term ``$G_n$-equivariant'' is reserved for $G_n$-vector bundles 
with algebraic $G_n$-action on their total spaces.) 
The composition of ${\mathcal S}$ with the generic fibre functor is 
the identical full embedding of $\mathfrak{SL}^u_n$ into the category  
of finite-dimensional $K_n$-semi-linear $G_n$-representations. 

In this section we show that the category $\mathfrak{SL}^u_n$ is 
abelian and its irreducible objects are generic fibres of irreducible 
coherent $G_n$-{\sl equivariant} sheaves on ${\mathbb P}^n_k$, i.e., 
direct summands of ${\rm Hom}_{K_n}((\Omega^n_{K_n/k})^{\otimes M},
\bigotimes_{K_n}^{\bullet}\Omega^1_{K_n/k})$ 
for appropriate integer $M\ge 0$. 

\begin{lemma} The category $\mathfrak{SL}^u_n$ is closed 
under taking $K_n$-semi-linear subquotients.\end{lemma} 
{\it Proof.} Let $V\in\mathfrak{SL}^u_n$ and $0\to V_1\to V
\stackrel{\pi}{\to}V_2\to 0$ be a short exact sequence of 
semi-linear representations of $G_n$ over $K_n$. As the $k$-vector space 
$V^{T_n^{{\rm tors}}}$ (of the elements in $V$ fixed by the torsion 
subgroup $T_n^{{\rm tors}}$ in $T_n$) spans the $K_n$-vector space $V$, 
the $k$-vector space $\pi(V^{T_n^{{\rm tors}}})\subseteq 
V_2^{T_n^{{\rm tors}}}$ spans the $K_n$-vector space $V_2$. 

This means that $V_2=V_2^{T_n^{{\rm tors}}}\otimes_kK_n$ and 
$\pi(V^{T_n^{{\rm tors}}})=V_2^{T_n^{{\rm tors}}}$. 

In other words, the sequence of $T_n^{{\rm tors}}$-invariants 
$0\to V_1^{T_n^{{\rm tors}}}\to V^{T_n^{{\rm tors}}}\to 
V_2^{T_n^{{\rm tors}}}\to 0$ is exact, and extending its coefficients 
to $K_n$ gives the exact sequence $0\to V_1^{T_n^{{\rm tors}}}\otimes_k
K_n\to V=V^{T_n^{{\rm tors}}}\otimes_kK_n\stackrel{\pi'}{\to}V_2=
V_2^{T_n^{{\rm tors}}}\otimes_kK_n\to 0$. As $\pi$ coincides with 
$\pi'$, we get $V_1=V_1^{T_n^{{\rm tors}}}\otimes_kK_n$. 

Clearly, any subquotient of a unipotent representation of $T_n$ 
is again unipotent, and thus, $V_1,V_2\in\mathfrak{SL}^u_n$. \qed 

\begin{lemma} \label{irre-equi} Let $E$ be the total space of a vector 
bundle on ${\mathbb P}^n_k$, ${\rm Aut}_{{\rm lin}}(E)$ be the group 
of automorphisms of $E$ over $k$ inducing linear transforms between 
the fibres, and $\tau:G_n\to{\rm Aut}_{{\rm lin}}(E)$ be an irreducible 
$G_n$-structure on $E$, i.e., a discrete group homomorphism 
splitting the projection ${\rm Aut}_{{\rm lin}}(E)\to G_n$. 
Then the Zariski closure $\overline{\tau(G_n)}$ is reductive. \end{lemma} 
{\it Proof.} Let ${\rm Aut}_{\tau}$ be the kernel of the projection 
$\overline{\tau(G_n)}\stackrel{\pi}{\to}G_n$. 

For each point $p\in{\mathbb P}^n_k$ let $\rho_p:R_p:=\pi^{-1}
({\rm Stab}_p)\to{\rm GL}(E_p)$ be the natural representation. 

As we suppose that $E$ is an irreducible $G_n$-bundle, $\rho_p$ is 
irreducible, since otherwise $B:=\overline{\tau(G_n)}B_p\subset E$ is a 
$G_n$-subbundle for any proper $R_p$-invariant $k$-subspace $B_p\subset E_p$. 

In particular, $\rho_p$ is trivial on the unipotent radical of 
$R_p$. The unipotent radical of any algebraic 
group contains the unipotent radical of its arbitrary normal subgroup, 
so $\rho_p$ is trivial on the unipotent radical of ${\rm Aut}_{\tau}$. 
As the action of ${\rm Aut}_{\tau}$ on $E$ is faithful, $\bigcap_p
\ker\rho_p|_{{\rm Aut}_{\tau}}=\{1\}$, i.e., ${\rm Aut}_{\tau}$ is 
reductive. As $G_n$ is also reductive, so is $\overline{\tau(G_n)}$. \qed 

\vspace{4mm}

For a commutative finite $k$-algebra $A$ denote by $R_{A/k}$ the Weil functor 
of restriction of scalars on $A$-schemes, cf. \cite{dg}, I, \S 1, 6.6. 

We need the following particular case of Th\'eor\`eme 8.16 of \cite{bt}.
\begin{theorem} \label{BT} Let $G$ be a simply connected absolutely 
almost simple $k$-group, and $G'$ be a reductive $k$-group. Let 
$\tau:G(k)\to G'(k)$ be a homomorphism with Zariski dense image. 
Let $G'_1,\dots,G'_m$ be the almost simple normal subgroups of $G'$. 

Then there exist finite field extensions $k_i/k$, field embeddings 
$\varphi_i:k\to k_i$, a special isogeny $\beta:\prod_{i=1}^mR_{k_i/k}
{}^{\varphi_i}G\to G'$ (here ${}^{\varphi_i}G:=G\times_{k,\varphi_i}k_i$) 
and a homomorphism $\mu:G(k)\to Z_{G'}(k)$ such that $\beta(R_{k_i/k}
{}^{\varphi_i}G)=G'_i$ and $\tau(h)=\mu(h)\cdot\beta(\prod_{i=1}^m
\varphi^{\circ}_i(h))$ for any $h\in G(k)$ (here $\varphi^{\circ}_i:G(k)
\to(R_{k_i/k}{}^{\varphi_i}G)(k)$ is the canonical homomorphism). \qed 
\end{theorem}

\begin{corollary} \label{cor-BT} Under assumptions of Theorem \ref{BT}, 
for any torus $T\subset G$ the Zariski closure of $\tau(T(k))$ is a torus 
in $G'$. \qed \end{corollary}

\begin{proposition} \label{equivar-irred} If $n\ge 2$ then any irreducible 
object of $\mathfrak{SL}^u_n$ is a direct summand of ${\rm Hom}_{K_n}
((\Omega^n_{K_n/k})^{\otimes M},\bigotimes_{K_n}^{\bullet}
\Omega^1_{K_n/k})$ for an appropriate $M$. \end{proposition}
{\it Proof.} The functor ${\mathcal S}$, mentioned in the beginning of 
this \S, associates to an irreducible object $V$ of $\mathfrak{SL}^u_n$ 
a coherent $G_n$-sheaf ${\mathcal V}$ on ${\mathbb P}^n_k$ with generic 
fibre $V$. 

Let, as before, $T_n$ be a maximal torus in $G_n$ and 
$Y_n\subset{\mathbb P}^n_k$ be the $n$-dimensional $T_n$-orbit. 
As $V^{T_n^{{\rm tors}}}=\Gamma(Y_n,{\mathcal V})^{T_n^{{\rm tors}}}$, cf. 
\cite{pgl}, is a unipotent representation of $T_n$, Lemma \ref{irre-equi} 
and Corollary \ref{cor-BT} imply that $\Gamma(Y_n,{\mathcal V}) 
^{T_n^{{\rm tors}}}$ is a trivial representation of $T_n$. 

In a $k$-basis of $V^{T_n^{{\rm tors}}}$ the $G_n$-action on $V$ determines 
a 1-cocycle $(g_{\sigma})\in Z^1(G_n,{\rm GL}_MK_n)$, where $M=\dim_{K_n}V$. 
There is an integer $N>n+2$ and elements $\alpha_1,\dots,\alpha_N\in G_n$ 
such that the morphism $(T_n)^N\stackrel{\pi}{\longrightarrow}G_n$, given 
by $(h_1,\dots,h_N)\mapsto\alpha_1h_1\alpha_1^{-1}\cdots\alpha_Nh_N
\alpha_N^{-1}$, is surjective. Namely, using the Gau\ss\ elimination 
algorithm, one shows that any element of $G_n$ is a product of $\le(n+1)^2$ 
elementary matrices and an element of $T_n$. On the other hand, it 
follows from the identity $\left(\begin{array}{rc} 1&0\\ -1&1\end{array}
\right)\left(\begin{array}{cc} a&0\\ 0&b\end{array}\right)\left(
\begin{array}{cc} 1&0\\ 1&1\end{array}\right)=\left(\begin{array}{cc} 
a&0\\ b-a&b\end{array}\right)$ that for any elementary matrix $\alpha$ 
the product $T_n\cdot\alpha T_n\alpha^{-1}\cdot T_n$ contains all 
elementary matrices of the same type as $\alpha$. This gives a 
surjection $(T_n\times\prod_{i\neq j}\alpha_{ij}T_n\alpha_{ij}^{-1})^{(n+1)^2}
\times T_n\stackrel{\times}{\longrightarrow\hspace{-3mm}
\to}G_n$, where $\alpha_{ij}$ is the elementary matrix 
with $1$ in the $i$-th row and $j$-th column. 

Then \begin{multline*} g_{\pi(h_1,\dots,h_N)}
=g_1(x)g'_1(\alpha_1h_1(x))g_2(\alpha_1h_1\alpha_1^{-1}(x))g'_2(\alpha_1h_1
\alpha_1^{-1}h_2(x))\cdots\\ g_N(\alpha_1h_1\alpha_1^{-1}\cdots
\alpha_{N-1}h_{N-1}\alpha_{N-1}^{-1}(x))g'_N(\alpha_1h_1\alpha_1^{-1}
\cdots\alpha_{N-1}h_{N-1}\alpha_{N-1}^{-1}\alpha_Nh_N(x)),\end{multline*} where 
$g_j(x):=g_{\alpha_j}$ and $g'_j(x):=g_{\alpha_j^{-1}}$ for all 
$1\le j\le N$. In other words, the lifting of the $G_n$-action on 
${\rm tot}({\mathcal V})$ to $(T_n)^N$-``coupling'' via $\pi$ determines 
a rational map $(T_n)^N\times{\mathbb P}^n_k\dasharrow{\rm GL}_Mk$. 
Clearly, it corresponds to a regular morphism $(T_n)^N\times
{\rm tot}({\mathcal V})\longrightarrow{\rm tot}({\mathcal V})$ and 
factors through a regular morphism $G_n\times{\rm tot}({\mathcal V})
\longrightarrow{\rm tot}({\mathcal V})$ of $k$-varieties, i.e., 
we see that ${\mathcal V}$ is equivariant. 

The generic fibres of irreducible $G_n$-equivariant sheaves 
on ${\mathbb P}^n_k$ are exactly of the desired type. \qed 

\vspace{4mm}

{\sc Remark.} There can exist, a priori, non-equivariant irreducible 
coherent $G_n$-sheaves on ${\mathbb P}^n_k$, e.g. the extension of 
coefficients to ${\mathcal O}_{{\mathbb P}^n_k}$ of a non-rational 
representation of $G_n$ is seemingly of this type. 

\section{``Positivity''} \label{polozhit}
In this section we show that for any admissible $F$-semi-linear 
representation $V$ of $G_{F/k}$ any irreducible subquotient of the 
$K_n$-semi-linear representation $V^{G_{F/K_n}}$ of $G$ is a direct 
summand of $\bigotimes_{K_n}^{\bullet}\Omega^1_{K_n/k}$. 

It is shown in \cite{pgl} that any finite-dimensional $K_n$-semi-linear 
$G_n$-representation extendable to ${\rm End}(K_n/k)$, e.g. $V^{G_{F/K_n}}$, 
is an object of $\mathfrak{SL}^u_n$. By Proposition \ref{equivar-irred}, 
we only need to eliminate the negative twists by $\Omega^n_{K_n/k}$ in 
irreducible subquotients of $V^{G_{F/K_n}}$. 

To do that we show first that the generic fibres of irreducible 
coherent $G_n$-equivariant sheaves are determined by their restrictions 
to the subgroup ${\rm Aff}_n=G_n\cap{\rm Aut}({\mathbb A}^n_k/k)$. 

\begin{lemma} \label{rest-t-af} Let ${\rm Aff}_n$ be 
the group of affine transformations of an affine space 
${\mathbb A}^n_k$ with the function field $K_n$. 
Then the natural morphism \begin{equation}\label{ext-coeff}
\{\mbox{{\rm rational $k$-linear ${\rm Aff}_n$-representations}}\}
\stackrel{\otimes_kK_n}{\longrightarrow}\{\mbox{{\rm $K_n$-semi-linear 
${\rm Aff}_n$-representations}}\}\end{equation} 
transforms isomorphism classes of irreducible $k$-representations 
of ${\rm Aff}_n$ to isomorphism classes of irreducible 
$K_n$-semi-linear representations of ${\rm Aff}^{(1)}_n{\mathbb Q}$, 
the subgroup of ${\rm Aff}_n$ consisting of ${\mathbb Q}$-affine 
substitutions of $x_1,\dots,x_n$ with Jacobian equal to 1. \end{lemma} 
{\it Proof.} Let $W$ be an irreducible $k$-representations of 
${\rm Aff}_n$, and $U\subset W\otimes_kK_n$ a non-zero 
$K_n$-semi-linear subrepresentation of ${\rm Aff}^{(1)}_n{\mathbb Q}$. 
Let $\alpha=\sum_{j=1}^Nw_j\alpha_j\in U$ be a non-zero element with 
minimal possible $N$, where $w_j\in W$ and $\alpha_j\in K_n$. Multiplying 
$\alpha$ by an element of $K_n$ we may assume that all $\alpha_j$ 
are polynomials: $\alpha=\sum_Iw'_Ix^I$. Since $W$ is irreducible, 
the elements of the unipotent radical $({\rm Aff}_n)_u$ of ${\rm Aff}_n$, 
i.e., $\sigma$ such that $\sigma z-z\in k$ for any linear function 
$z$ on ${\mathbb A}^n_k$, act trivially on $W$. 

Applying an appropriate composition of difference operators 
$\sigma-\tau$ for some $\sigma,\tau$ in the unipotent radical of 
${\rm Aff}^{(1)}_n{\mathbb Q}$ to $\alpha$, we can lower the degrees 
of the polynomials $\alpha_j$ and eventually get a non-zero element of 
$W$. As $W=W_0\otimes_{{\mathbb Q}}k$ for an irreducible representation 
$W_0$ of ${\rm Aff}^{(1)}_n{\mathbb Q}$, any non-zero element of $W$ 
generates $W\otimes_kK_n$, which means that $U=W\otimes_kK_n$. \qed 

\begin{corollary} \label{restr-to-aff} Let ${\rm Aff}_n$, 
$({\rm Aff}_n)_u$, ${\mathbb A}^n_k$ and $K_n$ be as in Lemma 
\ref{rest-t-af}. Then $V\mapsto V^{({\rm Aff}_n)_u}$ gives a 
natural bijection \begin{equation}\label{semi-to-lin} \left\{
\begin{array}{c} \mbox{{\rm isomorphism classes of}}\\ 
\mbox{{\rm irreducible $G_n$-subrepresentations in}}\\ 
\mbox{{\rm $\bigoplus_M{\rm Hom}_K((\Omega^n_{K_n/k})^{\otimes M},
\bigotimes_{K_n}^{\bullet}\Omega^1_{K_n/k})$}}\end{array}\right\}
\stackrel{\sim}{\longrightarrow}
\left\{\begin{array}{c}\mbox{{\rm isomorphism classes of}}\\ 
\mbox{{\rm irreducible rational $k$-linear}}\\ 
\mbox{{\rm ${\rm Aff}_n$-representations}}\end{array}\right\}
\end{equation} such that its composition with the morphism 
(\ref{ext-coeff}) is the inclusion map. \qed \end{corollary} 

\vspace{4mm}

Let $W$ be an $(n+1)$-dimensional $k$-vector space, $L\subset W$ be 
a one-dimensional subspace, and $H_{{\rm lin}}=\ker[{\rm GL}(W,L)\to
{\rm GL}(W/L)]\cong k^{\times}\ltimes{\rm Hom}(W/L,L)$ be the group 
preserving $L$ and inducing the identity automorphism of $W/L$. 

\begin{lemma}\label{H-inv} For any Young diagram $\lambda$ with no 
columns of height $\ge n+1$ one has $$(S^{\lambda}W^{\vee}\otimes_k
(\det W)^{\otimes s})^{H_{{\rm lin}}}=\left\{\begin{array}{ll} 
S^{\lambda}(W/L)^{\vee} & \mbox{{\rm if} $s=0$,} \\ 0 & 
\mbox{{\rm otherwise}}\end{array}\right.$$ \end{lemma}
{\it Proof.} Denote by $X=AF(W)\cong{\rm GL}(W)/R_u(B)$, the variety 
of complete affine flags in $W$. An affine flag is a filtration 
$W_{\bullet}=(0=W_0\subset W_1\subset W_2\subset\dots\subset W_{n+1}=W)$ 
with $\dim_kW_j=j$ and a collection of $l_j\in W_j/W_{j-1}-\{0\}$. 

Let $Y\cong{\rm GL}(W)/B$ be the variety of complete linear flags 
in $W$. Then the natural projection $X\stackrel{\pi}{\longrightarrow}Y$ 
is a principal $({\mathbb G}_m)^{n+1}$-bundle, and there is a 
decomposition $\pi_{\ast}{\mathcal O}_X=\bigoplus_{\mu}{\mathcal M}(\mu)$ 
into a direct sum of invertible sheaves on $Y$, where $\mu$ runs over 
the group ${\mathbb Z}^{n+1}$ of characters of $({\mathbb G}_m)^{n+1}$, 
so ${\mathcal O}(X)=\bigoplus_{\mu}\Gamma(Y,{\mathcal M}(\mu))$. 

Set $X^{\circ}=\{(V_{\bullet},l_{\bullet})\in X~|~V_n\cap L=0
\Leftarrow\!\!\!\!\Rightarrow``l_{n+1}\in L"\}$. Then reduction modulo 
$L$ defines a principal $L^{\oplus n}\times{\mathbb G}_m$-bundle 
$X^{\circ}\longrightarrow AF(W/L)$, $(l_1,\dots,l_{n+1})\mapsto
(l_1,\dots,l_n)$. 

Let $X^{\circ}\stackrel{``l_{n+1}"}{\longrightarrow}L-\{0\}$ 
be the natural $H$- (or ${\mathbb G}_m$-) equivariant map, 
and $\overline{l}_{n+1}$ be the composition of $``l_{n+1}"$ 
with a fixed isomorphism $L-\{0\}\cong{\mathbb G}_m$. 

Set $SH:=H_{{\rm lin}}\cap{\rm SL}(W)={\rm Hom}(W/L,L)$. Then 
${\mathcal O}(X)^{SH}={\mathcal O}(X)\cap{\mathcal O}(X^{\circ})^{SH}
={\mathcal O}(X)\cap{\mathcal O}(AF(W/L))[\overline{l}_{n+1},
\overline{l}_{n+1}^{-1}]$, so ${\mathcal O}(X)(\mu)^{SH}={\mathcal O}(X)
\cap{\mathcal O}(AF(W/L))(\mu')\overline{l}_{n+1}^{\mu_{n+1}}$, where 
$\mu'\in{\mathbb Z}^n$ is the restriction of $\mu$ to the first $n$ 
multiples of $({\mathbb G}_m)^{n+1}$. 

For any $\mu$ this is an irreducible 
representation of ${\rm GL}(W/L)$, and thus, ${\mathcal O}(X)(\mu)^{SH}
={\mathcal O}(AF(W/L))(\mu')\overline{l}_{n+1}^{\mu_{n+1}}$ if 
${\mathcal O}(X)(\mu)\neq 0$. 

As any irreducible representation of ${\rm SL}(W)$ coincides with 
${\mathcal O}(X)(\mu)$ for some $\mu$, this implies that 
$(S^{\lambda}W^{\vee})^{SH}=S^{\lambda}(W/L)^{\vee}$. \qed

\begin{theorem} \label{main} For any $F$-semi-linear 
$G_{F/k}$-representation $V\in{\mathcal A}$ any irreducible subquotient 
of the $K_n$-semi-linear $G_n$-representation $V^{G_{F/K_n}}$ is a direct 
summand of $\bigotimes_{K_n}^{\bullet}\Omega^1_{K_n/k}$. \end{theorem}
{\it Proof.} Let $W={\mathbb A}^{n+1}_k$ be the vector space with 
coordinates $x_1,\dots,x_{n+1}$, so $k(W)=K_{n+1}$. By Proposition 
\ref{equivar-irred} and Corollary \ref{restr-to-aff}, the restrictions 
to ${\rm Aff}_{n+1}={\rm Aff}(W)$ of irreducible subquotients of the 
$K_{n+1}$-semi-linear $G_{n+1}$-representation $V^{G_{F/K_{n+1}}}$ 
are of type $(S^{\lambda}W^{\vee}\otimes(\det W)^{\otimes s})
\otimes_kK_{n+1}$ for a Young diagram $\lambda$ with no columns 
of height $n+1$ and some integer $s$, where ${\rm Aff}_{n+1}$ 
acts on $W$ via its reductive quotient ${\rm GL}(W)$. 

Let $H\subset{\rm Aff}_{n+1}$ be the subgroup fixing the functionals 
$x_1,\dots,x_n$ in $W^{\vee}$ vanishing on $L$. Let ${\rm Aff}_u$ be the 
unipotent radical of ${\rm Aff}_{n+1}$, i.e. the group of translations of $W$. 

Then the restrictions to ${\rm Aff}_n$ of the irreducible subquotients of 
the $K_n$-semi-linear $G_n$-representation $V^{G_{F/K_n}}$ are contained 
in $((S^{\lambda}W^{\vee}\otimes(\det W)^{\otimes s})\otimes_kK_{n+1})^H$. 
As $H\cap{\rm Aff}_u=\langle 1\rangle_k\cong k$, we get $((S^{\lambda}
W^{\vee}\otimes(\det W)^{\otimes s})\otimes_kK_{n+1})^{H\cap{\rm Aff}_u}
=(S^{\lambda}W^{\vee}\otimes(\det W)^{\otimes s})\otimes_kK_n$, so 
$((S^{\lambda}W^{\vee}\otimes(\det W)^{\otimes s})\otimes_kK_{n+1})^H=
(S^{\lambda}W^{\vee}\otimes(\det W)^{\otimes s})^H\otimes_kK_n$. 

By Lemma \ref{H-inv}, $(S^{\lambda}W^{\vee}\otimes(\det W)^{\otimes s})^H$ 
coincides with $S^{\lambda}(W/L)^{\vee}$ if $s=0$, and vanishes otherwise. 
This means that any representation of ${\rm Aff}_n$ obtained this way is 
a direct summand of the tensor algebra of the representation $(W/L)^{\vee}
=(\Omega^1_{K_n/k})^{\{{\rm translations}\}}$ of ${\rm GL}_nk$. 
As any irreducible subquotient $U$ of the $K_n$-semi-linear 
$G_n$-representation $V^{G_{F/K_n}}$ is determined by its restriction 
$U|_{{\rm Aff}_n}$ to ${\rm Aff}_n$ and $U|_{{\rm Aff}_n}$ is a direct summand 
of $\bigotimes_{K_n}^{\bullet}\Omega^1_{K_n/k}$, the same holds for $U$. \qed 

\section{Extensions in $\mathfrak{SL}^u_n$ and in ${\mathcal A}$}
\label{rasshir-v-A} For an integer $\ell\ge 2$ such that $\mu_{\ell}
\subset k$ (see \S\ref{obozn}), denote by ${\rm Aff}_n^{(\ell)}{\mathbb Q}$ 
the subgroup of ${\rm Aff}_n$ consisting of the ${\mathbb Q}
(\mu_{\ell})$-affine substitutions of $x_1,\dots,x_n$ with Jacobian 
in $\mu_{\ell}$: $x_i\mapsto\sum_{j=1}^na_{ij}x_j+b_i$, where 
$a_{ij},b_i\in{\mathbb Q}(\mu_{\ell})\subset k$ and $\det(a_{ij})\in
\mu_{\ell}$; and by ${\rm SAff}_n^{(\ell)}{\mathbb Q}$ the subgroup 
of index $\ell$ consisting of elements with Jacobian equal to 1: 
$\det(a_{ij})=1$. 
\begin{lemma} \label{aff-spl} Let $n,\ell\ge 2$ be integers. Assume 
that $\mu_{\ell}\subset k$. Let $U_0$ be the unipotent radical of 
${\rm SAff}_n^{(\ell)}{\mathbb Q}$. Then for any object $V\in
\mathfrak{SL}^u_n$ there is a rational representation $W$ of the 
reductive quotient ${\rm SL}_n{\mathbb Q}(\mu_{\ell})=
{\rm SAff}_n^{(\ell)}{\mathbb Q}/U_0$ of ${\rm SAff}_n^{(\ell)}{\mathbb Q}$, 
and an isomorphism of semi-linear ${\rm SAff}_n^{(\ell)}{\mathbb Q}$-modules 
$W\otimes_{{\mathbb Q}(\mu_{\ell})}K_n\stackrel{\sim}{\longrightarrow}V$. 

Irreducible rational representations of ${\rm SL}_n{\mathbb Q}(\mu_{\ell})$ 
with coefficients extended to $K_n$ are irreducible 
semi-linear representations of ${\rm SAff}_n^{(\ell)}{\mathbb Q}$ 
over $K_n$. In particular, any extension 
in $\mathfrak{SL}^u_n$ splits as an extension of $K_n$-semi-linear 
representations of ${\rm SAff}_n^{(\ell)}{\mathbb Q}$. \end{lemma}
{\it Proof.} It is shown in Lemma 6.3 (1) of \cite{pgl} that $H^0(U_0,-)$ 
is a fibre functor on $\mathfrak{SL}^u_n$ independent of $\ell$, 
so $V=V^{U_0}\otimes_kK_n$, i.e, the restriction of $V$ to 
${\rm SAff}_n^{(\ell)}{\mathbb Q}$ is a $k$-linear representation 
$V^{U_0}$ of ${\rm SL}_n{\mathbb Q}(\mu_{\ell})$ with 
coefficients extended to $K_n$, for any $V\in\mathfrak{SL}^u_n$. 

As it follows from Proposition \ref{equivar-irred}, the irreducible 
subquotients $V_{\alpha}$ of $V$ restricted to ${\rm SAff}_n^{(\ell)}
{\mathbb Q}$ are of the form $W_{\alpha}\otimes_{{\mathbb Q}(\mu_{\ell})}
K_n$, where $W_{\alpha}$ are rational irreducible 
representations of ${\rm SL}_n{\mathbb Q}(\mu_{\ell})$. 
Then the irreducible subquotients of $V^{U_0}$ are 
$V_{\alpha}^{U_0}=W_{\alpha}\otimes_{{\mathbb Q}(\mu_{\ell})}k$, 
and $V_{\alpha}$ are irreducible semi-linear representations of 
${\rm SAff}_n^{(\ell)}{\mathbb Q}$ by Lemma \ref{rest-t-af}. 

If $V^{U_0}$ is not semi-simple then it admits a non-semi-simple 
subquotient $W$ of length 2. Let in Theorem \ref{LR} 
$\kappa={\mathbb Q}(\mu_{\ell})$, $K=k$, $G={\rm SL}_n$, 
${\mathcal G}$ be the Zariski closure of the image of 
${\rm SL}_n{\mathbb Q}(\mu_{\ell})$ in ${\rm GL}_k(W)$ and let 
$\tau$ be given by the ${\rm SL}_n{\mathbb Q}(\mu_{\ell})$-action 
on $W$. Then the unipotent radical of ${\mathcal G}$ is commutative. 
As the derivations of $\kappa={\mathbb Q}(\mu_{\ell})$ are zero, 
we see that the $k$-linear representation $W$ of 
${\rm SL}_n{\mathbb Q}(\mu_{\ell})$ is semi-simple. \qed 

\vspace{4mm}

{\sc Remarks.} 1. Using Theorems \ref{BT} and \ref{LR} it is not 
hard to show that any representation of ${\rm SL}_nK$ over any field 
of characteristic zero is semi-simple for any number field $K$. 

2. Let $V=\Omega^1_{K_n}/\Lambda\otimes_kK_n$, where 
$\Lambda\subset\Omega^1_k$ is a proper $k$-subspace. Let the extension 
$0\to V\to U\to K_n\to 0$ be given by the cocycle $(\omega_{\sigma}=
d\log\frac{\sigma\omega}{\omega})\in Z^1(G_n,V)$, where 
$\omega=dx_1\wedge\dots\wedge dx_n\in\Omega^n_{K_n/k}$. Then 
the restriction of $(\omega_{\sigma})$ to ${\rm GL}_nk$ is non-trivial. 

3. The convolution with the Euler vector field $\sum_{j=1}^nx_j
\partial/\partial x_j$ defines a ${\rm GL}_nk$-equivariant morphism 
$\Omega^1_{K_n/k}\to K_n$ given by $dx_j\mapsto x_j$. It is non-split 
for $n\ge 3$, since $(\Omega^1_{K_n/k})^{{\rm SL}_n{\mathbb Q}}=0$. 

\begin{lemma} \label{usl-dok-corr} Let $n,\ell\ge 2$ and $s$ be some 
integers, and $\lambda$ be a Young diagram with columns of height $<n$ 
such that $\ell$ does not divide $s+\frac{|\lambda|}{n-1}$, if $\lambda$ 
is rectangular of height $n-1$ and non-empty. Let $V=S^{\lambda}
_{K_n}\Omega^1_{K_n/k}\otimes_{K_n}(\Omega^n_{K_n/k})^{\otimes s}$. 

Then $(V^{{\mathcal H}^{(\ell)}_n})^{{\rm Aff}^{(\ell)}_{n-1}{\mathbb Q}}
=V^{{\rm Aff}^{(\ell)}_n{\mathbb Q}}$, where ${\mathcal H}^{(\ell)}_n
:=G_{K_n/K_{n-1}}\cap{\rm  Aff}^{(\ell)}_n{\mathbb Q}$. \end{lemma}
{\it Proof.} Let $W$ be the $k$-span of $dx_1,\dots,dx_n$ in 
$\Omega^1_{K_n/k}$. Then $S^{\lambda}_{K_n}\Omega^1_{K_n/k}=
S^{\lambda}_kW\otimes_kK_n$ and $\Omega^n_{K_n/k}=\det_kW\otimes_kK_n$. 
Set $S{\mathcal H}^{(\ell)}_n={\mathcal H}^{(\ell)}_n\cap
{\rm SAff}^{(\ell)}_n{\mathbb Q}$. Then ${\mathcal H}^{(\ell)}_n
\cong\mu_{\ell}\ltimes S{\mathcal H}^{(\ell)}_n$, and therefore, 
$V^{{\mathcal H}^{(\ell)}_n}=(V^{S{\mathcal H}^{(\ell)}_n})^{\mu_{\ell}}$. 

One has $V^{S{\mathcal H}^{(\ell)}_n}=(S^{\lambda}_kW\otimes_kK_n)
^{S{\mathcal H}^{(\ell)}_n}\otimes_k(\det_kW)^{\otimes s}$. As the 
intersection of the unipotent radical of ${\rm  Aff}^{(\ell)}_n
{\mathbb Q}$ with $S{\mathcal H}^{(\ell)}_n$ (i.e. the ${\mathbb Q}
(\mu_{\ell})$-translations of $x_n$) acts trivially on 
$S^{\lambda}_kW$ and fixes exactly $K_{n-1}$ in $K_n$, if $n\ge 1$, 
we get $$V^{S{\mathcal H}^{(\ell)}_n}=(S^{\lambda}_kW)
^{S{\mathcal H}^{(\ell)}_n}\otimes_k(\det_kW)^{\otimes s}\otimes_k
K_{n-1}=S^{\lambda}_k(W^{S{\mathcal H}^{(\ell)}_n})\otimes_k
(\det_kW)^{\otimes s}\otimes_kK_{n-1}.$$ Then 
\begin{multline*}V^{{\mathcal H}^{(\ell)}_n}=S^{\lambda}_k
(W^{S{\mathcal H}^{(\ell)}_n})\otimes_k
((\det_kW)^{\otimes s})^{\mu_{\ell}}\otimes_kK_{n-1}\\ 
=\left\{\begin{array}{ll}S^{\lambda}_{K_{n-1}}\Omega^1_{K_{n-1}/k}
\otimes_{K_{n-1}}(\Omega^{n-1}_{K_{n-1}/k})^{\otimes s} & 
\mbox{if $\ell|s$}\\ 0 & \mbox{otherwise.}\end{array}\right.\end{multline*} 

On the other hand, $V^{{\rm Aff}^{(\ell)}_n{\mathbb Q}}=
(S^{\lambda}_kW\otimes_k(\det_kW)^{\otimes s}\otimes_kK_n)
^{{\rm Aff}^{(\ell)}_n{\mathbb Q}}$ coincides with $(S^{\lambda}_kW
\otimes_k(\det_kW)^{\otimes s})^{{\rm Aff}^{(\ell)}_n{\mathbb Q}}$, 
since the unipotent radical of ${\rm Aff}^{(\ell)}_n{\mathbb Q}$ 
acts trivially on $S^{\lambda}_kW\otimes_k(\det_kW)^{\otimes s}$ 
and fixes exactly $k$ in $K_n$. Thus, for $n\ge 1$, we get 
$V^{{\rm Aff}^{(\ell)}_n{\mathbb Q}}=\left\{
\begin{array}{ll} k & \mbox{if $\lambda=0$ and $\ell|s$,}\\ 
0 & \mbox{otherwise.}\end{array}\right.$ This implies that 
$(V^{{\mathcal H}^{(\ell)}_n})^{{\rm Aff}^{(\ell)}_{n-1}{\mathbb Q}}=
\left\{\begin{array}{ll} k & \mbox{if $\lambda$ is rectangular 
of height $n-1$ and $\ell|(s+\frac{|\lambda|}{n-1})$,} \\ 
0 & \mbox{otherwise}\end{array}\right.$ for $n\ge 2$ 
(assuming that empty $\lambda$ is $(0\times(n-1))$-rectangular). \qed

\begin{lemma}[\cite{pgl}, Lemma 7.1] \label{quasi-desc} 
Let $n>m\ge 0$ be integers and $H$ be a subgroup of $G_{F/k}$ 
preserving $K_n$ and projecting onto a subgroup of $G_{K_n/k}$ 
containing the permutation group of the set $\{x_1,\dots,x_n\}$. 
Then the subgroup in $G_{F/k}$ generated by $G_{F/K_m}$ 
and $H$ is dense. \qed \end{lemma}

We note that ${\rm Aff}^{(\ell)}_n{\mathbb Q}\subset G_n
\subset G_{K_n/k}$ does indeed contain the permutation group 
of the set $\{x_1,\dots,x_n\}$ for any even $\ell\ge 2$. 

For any $U\in{\mathcal A}$ and $m\ge 0$ set $U_m=U^{G_{F/K_m}}$. Using 
smooth cochains, one defines the smooth cohomology $H^j_{{\rm smooth}}
(G_{F/k},-):={\rm Ext}^j_{{\mathcal S}m_{G_{F/k}}}({\mathbb Q},-)$. 
\begin{proposition} \label{proj-F} If $U\in{\mathcal A}$ and there is 
a subquotient of $U_n\in\mathfrak{SL}^u_n$ isomorphic to $K_n$ then 
there is an embedding $F\hookrightarrow U$ in ${\mathcal A}$. One 
has $H^1_{{\rm smooth}}(G_{F/k},V)=0$ for any $V\in{\mathcal A}$. 
\end{proposition}
{\it Proof.} By Lemma \ref{quasi-desc}, $U^{G_{F/k}}=
U_{n+1}^{{\rm Aff}^{(\ell)}_{n+1}{\mathbb Q}}\cap U_n$ 
for any even $\ell\ge 2$. By Theorem \ref{main}, Lemma 
\ref{aff-spl} and Lemma \ref{usl-dok-corr}, 
$(U_{n+1}^{{\mathcal H}^{(\ell)}_{n+1}})^{{\rm Aff}^{(\ell)}_n{\mathbb Q}}
=U_{n+1}^{{\rm Aff}^{(\ell)}_{n+1}{\mathbb Q}}$ for any $n\ge 1$ and any 
sufficiently big $\ell$ (where ${\mathcal H}^{(\ell)}_n$ is defined in 
Lemma \ref{usl-dok-corr}). Then, as $U_n\subseteq 
U_{n+1}^{{\mathcal H}^{(\ell)}_{n+1}}$, one has $U^{G_{F/k}}=
U_n^{{\rm Aff}^{(\ell)}_n{\mathbb Q}}$ for any sufficiently big even 
$\ell$, and thus, $U^{G_{F/k}}\neq 0$ if there is a subquotient of 
$U_n\in\mathfrak{SL}^u_n$ isomorphic to $K_n$.  

Clearly, ${\rm Ext}^j_{{\mathcal S}m_{G_{F/k}}}({\mathbb Q},-)=
{\rm Ext}^j_{{\mathcal C}}(F,-)$ on ${\mathcal C}$ for any 
$j\ge 0$,\footnote{Any class in ${\rm Ext}^j_{{\mathcal C}}(F,V)$ 
represented by $0\to V\to V_j\to\dots\to V_1\to F\to 0$ is sent to 
the class of $0\to V\to V_j\to\dots\to V_2\to V_1\times_F{\mathbb Q}
\to{\mathbb Q}\to 0$ in ${\rm Ext}^j_{{\mathcal S}m_{G_{F/k}}}
({\mathbb Q},V)$. Conversely, the class of $0\to V\to U_j\to\dots\to 
U_1\to{\mathbb Q}\to 0$ in ${\rm Ext}^j_{{\mathcal S}m_{G_{F/k}}}
({\mathbb Q},V)$ is sent to the class of $0\to V\to(U_j\otimes F)/K
\to U_{j-1}\otimes F\to\dots\to U_1\otimes F\to F\to 0$, where $K$ 
is the kernel of the surjection ${\rm forget}(V)\otimes F\to V$ and 
${\rm forget}:{\mathcal C}\longrightarrow{\mathcal S}m_{G_{F/k}}$ 
is the forgetful functor.} so we have to show that any smooth 
$F$-semi-linear extension $0\longrightarrow V\longrightarrow 
U\longrightarrow F\longrightarrow 0$ splits. 

Fix $u\in U$ in the preimage 
of $1\in F$. The stabilizer of $u$ contains a subgroup of type $G_{F/L}$ 
such that the elements of $L$ are algebraic over $K_m$ for some $m>1$. 
Then the normalized trace ${\rm tr}_{/K_m}u\in U_m$ belongs again to 
the preimage of $1\in K_m$, so $U_m$ surjects onto $K_m$. 

By Theorem \ref{main} and Lemma \ref{aff-spl} the semi-linear 
representation $U_m$ of ${\rm Aff}^{(\ell)}_m{\mathbb Q}$ over $K_m$ 
splits as $K_m\oplus V_m$, and thus, $U^{G_{F/k}}$ projects onto $k$. 
Then sending $1\in k\subset F$ to one of its preimages in 
$U^{G_{F/k}}$ extends to a splitting of $U\longrightarrow F$. \qed 

\begin{corollary} \label{exam-non-adm} For an integer $n\ge 1$ let 
$H\subseteq G_{F/k}$ be a subgroup containing $G_{F/K_n}$ such that 
$G_{F/\overline{K_n}}$ is a normal subgroup in $H$. Consider $H/G_{F/K_n}$ 
as a subset in the set $\{K_n\stackrel{/k}{\hookrightarrow}\overline{K_n}\}$ of field 
embeddings of $K_n$ into its algebraic closure in $F$ over $k$. 
Suppose that $H/G_{F/K_n}$ contains ${\rm Aff}_n$. Let 
$V=F[G_{F/k}/H]^{\circ}\in{\mathcal C}$ consist of formal 
degree-zero $F$-linear combinations of elements in 
$G_{F/k}/H$.\footnote{In particular, if $H=G_{\{F,\overline{K_n}\}/k}$ 
then $V\cong F[\{L\subset F~|~L\cong\overline{K_n}\}]^{\circ}$ 
consists of formal degree-zero $F$-linear combinations of algebraically 
closed subfields in $F$ of transcendence degree $n$ over $k$.} 
Then any quotient of $V$ which lies in ${\mathcal A}$ is zero. 
\end{corollary}
{\it Proof.} $V$ is generated by $\alpha=[1]-[\sigma]\in 
V_{2n}^{\langle({\rm Aff}_{2n})_u,T_{2n}\rangle}$, where $\sigma$ 
sends $x_j$ to $x_{2n+1-j}$ for each $1\le j\le 2n$. Any admissible 
semi-linear quotient of $V$ is generated by the image of $\alpha$, 
which is, by Propositions \ref{equivar-irred} and \ref{proj-F}, 
fixed by the whole $G_{F/k}$. On the other hand, $\sigma\alpha=-\alpha$, 
so any admissible semi-linear quotient of $V$ is zero. \qed 

\begin{theorem} \label{final} The category ${\mathcal A}$ is abelian. 
The functor $H^0(G_{F/L},-)$ is exact on ${\mathcal A}$ 
for any subfield $L$ in $F$ containing $k$. \end{theorem}
{\it Proof.} We have to check that ${\mathcal A}$ is stable under 
taking quotients. Let $V\in{\mathcal A}$ and $V\stackrel{\pi}{\to}V'$ 
be a surjection of $F$-semi-linear representations of $G_{F/k}$. By 
Proposition \ref{proj-F}, for any $K\subset F$ of finite type over $k$ 
and any $v\in(V')^{G_{F/\overline{K}}}-\{0\}$, the extension 
$0\to\ker\pi\to\pi^{-1}(F\cdot v)\to F\to 0$ of $F$-semi-linear 
representations of $G_{F/\overline{K}}$ splits. This implies that the 
natural projection $V^{G_{F/K}}\stackrel{\pi_K}{\longrightarrow}
(V')^{G_{F/K}}$ is surjective, and thus, $V'$ is also an admissible 
semi-linear representation. 

The functor $H^0(G_{F/L},-)$ on ${\mathcal A}$ is the composition 
of the forgetful functor $\Phi:{\mathcal A}_k\to{\mathcal C}_L$, 
the functor $H^0(G_{F/\overline{L}},-)$ on ${\mathcal C}_L$ 
and the exact functor $H^0(G_{\overline{L}/L},-)$ on 
${\mathcal S}m_{G_{\overline{L}/L}}$. If $L$ is of finite transcendence 
degree over $k$ then the forgetful functor $\Phi$ factors through 
${\mathcal A}_{\overline{L}}$, so the composition 
$H^0(G_{F/\overline{L}},-)\circ\Phi$ is exact. If $L$ is of infinite 
transcendence degree over $k$ then $H^0(G_{F/\overline{L}},-)$ 
induces an equivalence of categories ${\mathcal S}m_{G_{F/k}}
\stackrel{\sim}{\longrightarrow}{\mathcal S}m_{G_{\overline{L}/k}}$, 
so $H^0(G_{F/\overline{L}},-)$ is also exact. \qed 

\begin{corollary} $H^1_{{\rm smooth}}(G_{F/k},\Omega^{\bullet}
_{F/k,{\rm closed}})=H^1_{{\rm smooth}}(G_{F/k},
\Omega^{\bullet}_{F/k,{\rm exact}})=0$. \end{corollary}
{\it Proof.} By Proposition \ref{proj-F}, $H^1_{{\rm smooth}}(G_{F/k},
\Omega^{\bullet}_{F/k})=0$. Then a piece of the long cohomological 
sequence of the short exact sequence $0\to\Omega^q_{F/k,{\rm closed}}\to
\Omega^q_{F/k}\stackrel{d}{\to}\Omega^{q+1}_{F/k,{\rm exact}}\to 0$ looks 
as $H^0(G_{F/k},\Omega^{q+1}_{F/k,{\rm exact}})\to H^1_{{\rm smooth}}
(G_{F/k},\Omega^q_{F/k,{\rm closed}})\to H^1_{{\rm smooth}}(G_{F/k},
\Omega^q_{F/k})=0$. Evidently, $H^0(G_{F/k},\Omega^{q+1}_{F/k,{\rm exact}})
=0$, so $H^1_{{\rm smooth}}(G_{F/k},\Omega^{\bullet}_{F/k,{\rm closed}})=0$. 

Clearly, $H^0(G_{F/k},H^q_{{\rm dR}/k}(F))=0$.\footnote{Let $\omega\in
\Omega^q_{A/k}\subset\Omega^q_{F/k}$ represent a $G_{F/k}$-fixed 
element for a smooth finitely generated $k$-subalgebra $A\subset F$. 
Fix $\sigma\in G_{F/k}$ such that $A$ and $\sigma(A)$ are algebraically 
independent over $k$. Then $\omega-\sigma\omega=d\eta$ for some
$\eta\in\Omega^{q-1}_{B/k}$, where $B\subset F$ is a smooth finitely 
generated $(A\otimes_k\sigma(A))$-subalgebra. Fix a $k$-algebra 
homomorphism $\varphi:\sigma(A)\longrightarrow\overline{k}\subset F$ 
and extend $id\cdot\varphi:A\otimes_k\sigma(A)\longrightarrow 
A\otimes_k\overline{k}\subset F$ to $\psi:B\longrightarrow F$. Then 
$\psi$ induces a morphism of differential graded $k$-algebras 
$\psi_{\ast}:\Omega^{\bullet}_{B/k}\longrightarrow
\Omega^{\bullet}_{F/k}$ identical on $\Omega^{\bullet}_{A/k}$, so 
$\omega=d\psi_{\ast}(\eta)$.} A piece of the long cohomological 
sequence of short exact sequence $0\to\Omega^q_{F/k,{\rm exact}}\to
\Omega^q_{F/k,{\rm closed}}\to H^q_{{\rm dR}/k}(F)\to 0$ looks as 
$$H^0(G_{F/k},H^q_{{\rm dR}/k}(F))\longrightarrow 
H^1_{{\rm smooth}}(G_{F/k},\Omega^q_{F/k,{\rm exact}})\longrightarrow 
H^1_{{\rm smooth}}(G_{F/k},\Omega^q_{F/k,{\rm closed}})=0,$$ so 
$H^1_{{\rm smooth}}(G_{F/k},\Omega^{\bullet}_{F/k,{\rm exact}})=0$. \qed 

\subsection{Extensions in $\mathfrak{SL}^u_n$} 
Now we need the following particular case of Bott's theorem. 
\begin{theorem}[\cite{bott}, cf. also \cite{dema}]\label{bott} 
If ${\mathcal V}$ is an irreducible $G_n$-equivariant coherent 
sheaf on ${\mathbb P}^n_k$ then there exists at most one 
$j\ge 0$ such that $H^j({\mathbb P}^n_k,{\mathcal V})\neq 0$. 
If $H^j({\mathbb P}^n_k,{\mathcal V})^{G_n}\neq 0$ then 
${\mathcal V}\cong\Omega^j_{{\mathbb P}^n_k/k}$. \qed \end{theorem}

We also need the following explicit description of the homomorphisms 
in the case of commutative unipotent radicals of the target groups. 
It confirms general expectations, sketched in Remark 8.19 
of \cite{bt} and in \cite{t}, \S 5.1. 
\begin{theorem}[\cite{lr}, Theorem 3] \label{LR} Let $G$ be a simple 
simply connected Chevalley group over a field $\kappa$ of characteristic 
zero. Let ${\mathcal G}$ be a connected algebraic group over a field 
extension $K$ of $\kappa$. Let $\tau:G(\kappa)\to{\mathcal G}(K)$ be 
a homomorphism with Zariski dense image. Assume that the unipotent 
radical ${\mathcal G}_u$ of ${\mathcal G}$ is commutative and the 
composition $G(\kappa)\stackrel{\tau}{\to}{\mathcal G}(K)\to G'(K)$, 
where $G'={\mathcal G}/{\mathcal G}_u$, is induced by a rational 
$K$-homomorphism $\lambda:G\times_{\kappa}K\longrightarrow G'$. 

Then ${\mathcal G}_u$ splits over a finite field extension $L/K$ into 
a direct sum of $r$ copies of the adjoint representation of $G'$, so 
$r=\dim{\mathcal G}_u/\dim G'$. 

Let $A=\kappa[\varepsilon_1,\dots,\varepsilon_r]/(\varepsilon_1^2,\dots,
\varepsilon_r^2)$ and ${\mathcal H}=R_{A/\kappa}(G\times_{\kappa}A)\cong 
G\ltimes{\mathfrak g}^{\oplus r}$, where ${\mathfrak g}={\rm Lie}(G)$ is 
the adjoint representation of $G$. 

Then there exist derivations $\delta_1,\dots,\delta_r:\kappa\to L$ 
and an $L$-isogeny $\mu:{\mathcal H}\times_{\kappa}L\longrightarrow
{\mathcal G}\times_{\kappa}L$ such that $\tau=\mu\circ\eta_{\delta}$, 
where $\eta_{\delta}:G(\kappa)\longrightarrow{\mathcal H}(L)$ is 
induced by the ring embedding $id+\sum_{j=1}^r\varepsilon_j\delta_j:
\kappa\to A\otimes_{\kappa}L$. \qed \end{theorem} 

\begin{lemma} \label{ext-1-om} Let $n\ge 2$. Suppose that 
${\rm Ext}^1_{\mathfrak{SL}^u_n}(K_n,V_{\circ})\neq 0$ for some 
irreducible object $V_{\circ}$ of $\mathfrak{SL}^u_n$. Then either 
$V_{\circ}\cong\Omega^1_{K_n/k}$, or $V_{\circ}\cong{\rm Der}(K_n/k)$.  
One has ${\rm Ext}^1_{\mathfrak{SL}^u_n}(K_n,\Omega^1_{K_n/k})=k$ and 
${\rm Ext}^1_{\mathfrak{SL}^u_n}(K_n,{\rm Der}(K_n/k))={\rm Der}(k)$. 
\end{lemma}
{\it Proof.} Let ${\mathcal V}={\mathcal S}(V_{\circ})$ be the 
irreducible coherent $G_n$-equivariant sheaf on ${\mathbb P}(Q)
={\mathbb P}^n_k$ with the generic fibre $V_{\circ}$, and let 
$0\longrightarrow V_{\circ}\longrightarrow V\longrightarrow 
K_n\longrightarrow 0$ be an extension in $\mathfrak{SL}^u_n$. 

Suppose that the short exact sequence $0\to{\mathcal V}\to
{\mathcal S}(V)\to{\mathcal O}\to 0$ of coherent sheaves on 
${\mathbb P}(Q)$ splits. Let $E$ be the total space of 
${\mathcal S}(V)\cong{\mathcal O}\oplus{\mathcal V}$. 
Then, as ${\rm Aut}({\mathcal S}(V),{\mathcal V})\cong
({\mathbb G}_m\times{\mathbb G}_m)\ltimes
\Gamma({\mathbb P}^n_k,{\mathcal V})$, the $G_n$-structure 
on $V$ corresponds to a splitting of the sequence 
\begin{equation} \label{o-omega-2} 
1\longrightarrow({\mathbb G}_m\times{\mathbb G}_m)\ltimes 
\Gamma({\mathbb P}(Q),{\mathcal V})\longrightarrow{\rm Aut}_{{\rm lin}}
(E,{\rm tot}({\mathcal V}))\longrightarrow G_n\longrightarrow 1. 
\end{equation} 
As $H^1(G_n,{\mathbb G}_m\times{\mathbb G}_m)=1$, Theorem \ref{LR} 
(with $G={\rm SL}_{n+1}k$, $G'=G_n$ and ${\mathcal G}_u\subseteq
\Gamma({\mathbb P}^n_k,{\mathcal V})$) implies that a non-standard 
splitting of (\ref{o-omega-2}) can exist only if $\Gamma
({\mathbb P}^n_k,{\mathcal V})$ is isomorphic to the adjoint 
representation of $G_n$, i.e., if ${\mathcal V}\cong
{\mathcal T}_{{\mathbb P}^n_k/k}$. The identity 
${\rm Ext}^1_{\mathfrak{SL}^u_n}(K_n,{\rm Der}(K_n/k))
={\rm Der}(k)$ follows also from Theorem \ref{LR}.

If ${\mathcal V}\cong\Omega^1_{{\mathbb P}^n_k/k}$ then the target of 
the homomorphism ${\rm Ext}^1_{\mathfrak{SL}^u_n}
(K_n,{\mathcal V}_{\circ})\stackrel{\alpha}{\longrightarrow}
{\rm Ext}^1_{\mathcal{O}}(\mathcal{O}_{{\mathbb P}^n_k},{\mathcal V})=k$ 
induced by the functor ${\mathcal S}$ is generated by the 
class of the Euler extension $0\to\Omega^1_{{\mathbb P}(Q)/k}
\to Q^{\vee}\otimes_k{\mathcal O}_{{\mathbb P}(Q)}(-1)
\to{\mathcal O}_{{\mathbb P}(Q)}\to 0$. 
Let $E$ be the total space of the vector bundle with the sheaf of sections 
$Q^{\vee}\otimes_k{\mathcal O}(-1)$. Any $G_n$-structure on the middle term of 
this extension corresponding to an element of ${\rm Ext}^1_{\mathfrak{SL}^u_n}
(K_n,\Omega^1_{K_n/k})$ is a splitting of the short exact sequence 
\begin{equation} \label{o-omega} 1\longrightarrow{\rm Aut}(Q^{\vee}
\otimes_k{\mathcal O}(-1),\Omega^1_{{\mathbb P}(Q)/k})\longrightarrow
{\rm Aut}_{{\rm lin}}(E,{\rm tot}(\Omega^1_{{\mathbb P}(Q)/k}))
\longrightarrow G_n\longrightarrow 1. \end{equation} 
As any point of ${\mathbb P}(Q)$ determines a hyperplane in 
$Q^{\vee}$, the group ${\rm Aut}(Q^{\vee}\otimes_k{\mathcal O}(-1),
\Omega^1_{{\mathbb P}(Q)})$ coincides with the subgroup of 
${\rm GL}(Q)$ stabilizing all hyperplanes in $Q^{\vee}$, i.e., 
with the centre ${\mathbb G}_m$ of ${\rm GL}(Q)$. Then 
${\rm Aut}_{{\rm lin}}(E,{\rm tot}(\Omega^1_{{\mathbb P}(Q)/k}))$ 
is a central ${\mathbb G}_m$-extension of $G_n$, so the splitting 
of (\ref{o-omega}) is unique and corresponds to the usual 
$G_n$-equivariant structure. \qed

\begin{corollary} \label{number-case} If $k=\overline{{\mathbb Q}}$ 
then $\mathfrak{SL}^u_n$ is equivalent to the category of
$G_n$-equivariant vector bundles on ${\mathbb P}^n_k$. \end{corollary}
{\it Proof.} The category of $G_n$-equivariant vector bundles on
${\mathbb P}^n_k$ is a full sub-category of $\mathfrak{SL}^u_n$
with the same irreducible objects. As it is mentioned at the beginning 
of \S\ref{no-1}, p.\pageref{no-1}, the objects of $\mathfrak{SL}^u_n$ are 
generic fibres of coherent $G_n$-sheaves on ${\mathbb P}^n_k$. Suppose that 
$V\in\mathfrak{SL}^u_n$ is the generic fibre of a non-equivariant vector 
bundle on ${\mathbb P}^n_k$ of minimal possible rank. Then it fits into 
an exact sequence $0\to B\to V\to A\to 0$, where $A,B$ are the generic 
fibres of $G_n$-equivariant vector bundles on ${\mathbb P}^n_k$ 
and $A$ is irreducible. Let $0\neq C\subseteq B$ be an irreducible 
sub-object and $D=B/C$. Then the rows in the following commutative 
diagram are exact:
$$\begin{array}{ccccccccc} {\rm Hom}_u(A,D) &
\to & {\rm Ext}^1_u(A,C) & \to & {\rm Ext}^1_u(A,B) &
\to & {\rm Ext}^1_u(A,D) & \to & {\rm Ext}^2_u(A,C)\\
\| && \|\lefteqn{{\rm Lemma\ \ref{ext-1-om}}} && \bigcup && 
\|\lefteqn{{\rm minimality\ of}\ V} && \phantom{\xi}\uparrow\xi \\
{\rm Hom}_{{\rm eq}}(A,D) & \to & {\rm Ext}^1_{{\rm eq}}(A,C) 
& \to & {\rm Ext}^1_{{\rm eq}}(A,B) & \to & 
{\rm Ext}^1_{{\rm eq}}(A,D) & \to &
{\rm Ext}^2_{{\rm eq}}(A,C)\end{array}$$
where subscript $u$ refers to the category $\mathfrak{SL}^u_n$ 
and ${\rm eq}$ refers to the category of equivariant vector bundles. 

Let us show that $\xi$ is injective for any irreducible $A$ and $C$. 
As $C\otimes_{K_n}A^{\vee}$ is semi-simple (as it follows from 
Proposition \ref{equivar-irred}) and ${\rm Ext}^2_?(A,C)=
{\rm Ext}^2_?(K_n,C\otimes_{K_n}A^{\vee})$, where $?=u$ or ${\rm eq}$, 
we may assume that $A=K_n$ and $C$ is still irreducible. By Bott's 
Theorem \ref{bott}, $\dim_k{\rm Ext}^2_{{\rm eq}}(K_n,C)\le 1$ 
with equality only if $C\cong\Omega^2_{K_n/k}$, so we 
assume further that $C=\Omega^2_{K_n/k}$. The forgetful functor 
from $\mathfrak{SL}^u_n$ to the category of coherent sheaves on 
${\mathbb P}^n_k$ induces a homomorphism ${\rm Ext}^2_u(A,C)
\to H^2({\mathbb P}^n_k,\Omega^2_{{\mathbb P}^n_k/k})$. 
Clearly, its composition with $\xi$ is an isomorphism. 

Then the 5-lemma implies that ${\rm Ext}^1_u(A,B)=
{\rm Ext}^1_{{\rm eq}}(A,B)$, and thus, $V$ is equivariant. \qed 

\section{The category ${\mathcal A}$ in the case 
$k=\overline{{\mathbb Q}}$} \label{str-ra-A-Q}
In this section we determine (in Theorem \ref{irr-form}) the structure 
of the objects of ${\mathcal A}$ in the case $k=\overline{{\mathbb Q}}$, 
the field of algebraic numbers. The objects $V$ of ${\mathcal A}$ are 
quotients of sums of representations of $G$ over $k$ induced by rational 
representations of ${\rm GL}_mk$'s (considered as subquotients of $G$) 
with coefficients extended to $F$ (cf. Lemma \ref{semisimplicity}). 
Then we find (in Lemma \ref{first-rel}) a supply of elements in the 
induced representations vanishing in $V$, and use them in Lemmas 
\ref{calc}--\ref{to-tenz} to show that the objects of ${\mathcal A}$ 
are sums of quotients of $\bigotimes^{\bullet}_F{\mathfrak m}$. 
In \S\ref{proekt-pro-obra} we study extensions in ${\mathcal A}$. 

\vspace{4mm}

Let $V\in{\mathcal A}$ and $m\ge 0$ be such that $V_m\neq 0$. 
Then there is a non-zero morphism $F[G_{F/k}/G_{F/K_m}]
\otimes_{K_m[{\rm PGL}_{m+1}k]}V_m\longrightarrow V$ in ${\mathcal C}$. 
The object $V_m$ of $\mathfrak{SL}_m^u$ admits an irreducible sub-object 
$A\neq 0$. By Theorem \ref{main}, $A\cong S^{\lambda}_{K_m}\Omega^1_{K_m/k}$ 
for a Young diagram $\lambda$. Then $F[G_{F/k}/G_{F/K_m}]\otimes
_{K_m[{\rm PGL}_{m+1}k]}A\longrightarrow V$ is also non-zero. Clearly, 
$A=B\otimes_kK_m$, where $B:=A^{({\rm Aff}_m)_u}\cong(S^{\lambda}_{K_m}
\Omega^1_{K_m/k})^{({\rm Aff}_m)_u}\cong S^{\lambda}_k(k^m)$ is 
a rational irreducible representation of ${\rm GL}_mk:={\rm Aff}_m
/({\rm Aff}_m)_u$. 

This implies that there is a non-zero morphism $U:=
F[W^{\circ}]\otimes_{k[{\rm GL}_mk]}B\stackrel{\varphi}
{\longrightarrow}V$ and a surjection $U\longrightarrow 
S^{\lambda}_F\Omega^1_{F/k}$, where $W^{\circ}:=\{K_m\stackrel{/k}
{\hookrightarrow}F\}/({\rm Aff}_m)_u$ is considered as a $G_{F/k}$-set. 

As any embedding $K_m\stackrel{/k}{\hookrightarrow}F$ is determined 
by the images of $x_1,\dots,x_m$, one can consider $W^{\circ}$ as a 
subset of $(F/k)^m$ consisting of $m$-tuples with entries 
algebraically independent over $k$. 
More invariantly, let $W:={\rm Hom}_k((K_m/k)^{({\rm Aff}_m)_u},F/k)
\cong(F/k)^m$ be the group (a $k$-vector space) generated by $W^{\circ}$. 
The isomorphism is given by restriction of the homomorphisms to 
the basis $\{\overline{x_1},\dots,\overline{x_m}\}$ of 
$(K_m/k)^{({\rm Aff}_m)_u}$. 
Define a homogeneous map $\varkappa:W\longrightarrow\Omega^m_{F/k}
\otimes_k\det_k{\rm Hom}_G(F/k,W)$ of degree $m$ by inverting the 
first isomorphism in the sequence $W\stackrel{\sim}{\longleftarrow}
(F/k)\otimes_k{\rm Hom}_G(F/k,W)\stackrel{d\otimes id}{\hookrightarrow}
\Omega^1_{F/k}\otimes_k{\rm Hom}_G(F/k,W)\longrightarrow{\rm Sym}^m_F
(\Omega^1_{F/k}\otimes_k{\rm Hom}_G(F/k,W))\longrightarrow\Omega^m_{F/k}
\otimes_k\det_k{\rm Hom}_G(F/k,W)$. Then $W^{\circ}=\{w\in W~|
~\varkappa(w)\neq 0\}$. 

Let $(y_1,\dots,y_m)\mapsto[y_1,\dots,y_m]$ be the map 
$(F/k)^m\longrightarrow\{0\}\cup W^{\circ}$ sending $(y_1,\dots,y_m)$ 
to $[x_j\mapsto y_j]$ if $y_1,\dots,y_m$ are algebraically independent 
over $k$, and to 0 otherwise. Then $[\mu y_1,\dots,\mu y_m]\otimes b=
\mu^{|\lambda|}[y_1,\dots,y_m]\otimes b$ in $U$ for any $\mu\in k$. If 
$y_1,\dots,y_m$ belong to the $k$-linear envelope of $x_1,\dots,x_M$ 
for some integer $M\ge 1$ then $[y_1,\dots,y_m]\otimes b\in 
U_M^{({\rm Aff}_M)_u}$ is a weight $|\lambda|$ eigenvector 
of the centre of ${\rm GL}_Mk$. 

Let $U^!_M\subseteq U_M^{({\rm Aff}_M)_u}$ be the 
$k[{\rm GL}_Mk]$-envelope of $[x_1,\dots,x_m]\otimes b$ for some 
$b\neq 0$ (which is the same as $k$-envelope of all $[y_1,\dots,y_m]
\otimes c$ for algebraically independent $y_1,\dots,y_m$ in the $k$-linear 
envelope of $x_1,\dots,x_M$ and all $c$). Clearly, $U^!_m\cong B$ as 
$k[{\rm GL}_mk]$-modules, and any non-zero morphism $U\longrightarrow 
V$ induces an embedding $U^!_m\hookrightarrow V$.

\begin{lemma} \label{semisimplicity} If $k=\overline{{\mathbb Q}}$ 
then for any $V\in\mathfrak{SL}^u_M$ the representation 
$V^{({\rm Aff}_M)_u}$ of ${\rm GL}_Mk$ is rational semi-simple, 
and $V^{({\rm Aff}_M)_u}\otimes_kK_M=V$. \end{lemma}
{\it Proof.} By Corollary \ref{number-case}, $V$ is the generic 
fibre of a ${\rm PGL}_{M+1}k$-equivariant vector bundle. 

Then, by Lemma 6.3 (1) of \cite{pgl}, $V=V^{U_0}\otimes_kK_M$, 
where $U_0$ is a ${\mathbb Q}$-lattice in $({\rm Aff}_M)_u$. 
The group $({\rm Aff}_M)_u$ acts rationally on $V^{U_0}$. 
As the action of the ${\mathbb Q}$-lattice $U_0$ is trivial, 
the action of the entire $({\rm Aff}_M)_u$ is trivial, i.e., 
$V=V^{({\rm Aff}_M)_u}\otimes_kK_M$. Then the action of 
${\rm GL}_Mk$ on $V^{({\rm Aff}_M)_u}$ is rational, and thus, 
semi-simple. \qed 

\vspace{4mm}

{\sc Remark.} There is no semisimplicity 
if $k$ contains a transcendental element. 

Indeed, let the underlying $K_n$-vector 
space of $V\in\mathfrak{SL}^u_n$ be $K_n\oplus\Omega^1_{K_n}/
(\Lambda\otimes_kK_n)$ for a proper $k$-vector subspace 
$\Lambda\subset\Omega^1_k$ of finite codimension, and the 
$G_n$-action be given by $\sigma(f,\omega)=
(\sigma f,\sigma\omega+\sigma f\cdot d\log(\sigma\eta/\eta))$ for 
any $\sigma\in G_n$, where $\eta=dx_1\wedge\dots\wedge dx_n\in
\Omega^n_{K_n/k}$. (So $V$ fits into a non-split exact sequence 
$0\to\Omega^1_{K_n}/(\Lambda\otimes_kK_n)\to V\to K_n\to 0$.) 
Then $\sigma(f,\omega)=(\sigma f,\sigma\omega)$ if $\sigma
\in({\rm Aff}_n)_u$, and therefore, $V^{({\rm Aff}_n)_u}
=k\oplus\Omega^1_k/\Lambda$ is a non-trivial 
extension of trivial representations of ${\rm GL}_nk$. 

\begin{lemma} \label{first-rel} The kernel of $U^!_M
\longrightarrow S^{\lambda}_F\Omega^1_{F/k}$ is contained in the 
kernel of $U\stackrel{\varphi}{\longrightarrow}V$. \end{lemma}
{\it Proof.} By Lemma \ref{semisimplicity}, the image $\overline{U^!_M}$ 
of $U^!_M$ in $V$ is isomorphic to $\bigoplus_{|\nu|=|\lambda|}
(S^{\nu}_k(k^M))^{m_{\nu}}$. As $\overline{U^!_m}
\subseteq\overline{U^!_M}^{{\rm GL}_Mk\cap G_{K_M/K_m}}\cong
\bigoplus_{|\nu|=|\lambda|}(S^{\nu}_k(k^m))^{m_{\nu}}$ and 
$U^!_M$ is generated by $U^!_m\stackrel{\sim}{\longrightarrow}
S^{\lambda}_k(k^m)\subseteq S^{\lambda}_k(k^M)$, we see that 
$\overline{U^!_M}$ is isomorphic to $S^{\lambda}_k(k^M)$. 
Then the kernel of $U^!_M\longrightarrow S^{\lambda}_F
\Omega^1_{F/k}$ is contained in the kernel of the morphism 
$U\stackrel{\varphi}{\longrightarrow}V$. \qed

\vspace{5mm}

Let $W^{M\circ}\subset W^M$ be the subset consisting of $M$-tuples 
$(y_1,\dots,y_M)$ such that $\sum_{i\in I}y_i\in W^{\circ}$ for any 
non-empty subset $I\subseteq\{1,\dots,M\}$. 

Let $k[W^{M\circ}]\longrightarrow k[W^{\circ}]\otimes_{k[k^{\times}]}
k(M)$ be the $k$-linear map sending $(y_1,\dots,y_M)$ to $$\langle 
y_1,\dots,y_M\rangle:=\sum_{I\subseteq\{1,\dots,M\}}(-1)^{\# I}
[\sum_{i\in I}y_i]\in k[W^{\circ}]\otimes_{k[k^{\times}]}k(M).$$
Here $k(M)$ denotes a one-dimensional $k$-vector space with 
$k^{\times}$-action by $M$-th powers. As $(y,\dots,y)$ is sent 
to $$\sum_{j\ge 0}(-1)^j\binom{M}{j}j^M[y]=(t\frac{d}{dt})^M
(1-t)^M\vert_{t=1}\cdot[y]=(-1)^MM!\cdot[y],$$ 
it is surjective. Clearly, $\langle y_1,\dots,y_M\rangle
=\langle y_{\theta(1)},\dots,y_{\theta(M)}\rangle$ 
for any permutation $\theta\in\mathfrak{S}_M$. Let 
$\tilde{U}:=F[W^{|\lambda|\circ}]\longrightarrow U$ be the 
$F$-linear surjection sending $(y_1,\dots,y_{|\lambda|})$
to $\langle y_1,\dots,y_{|\lambda|}\rangle\otimes b$.

\begin{lemma} \label{calc} Let the $k$-linear map $\alpha:k[W^{\circ}]
\longrightarrow\bigotimes_k^MW$ be given by $[w]\mapsto w^{\otimes M}$. 
Then $\alpha$ factors through $k[W^{\circ}]\otimes_{k[k^{\times}]}k(M)$ 
and $\langle y_1,\dots,y_M\rangle\mapsto(-1)^M
\sum_{\theta\in\mathfrak{S}_M}y_{\theta(1)}\otimes\dots\otimes 
y_{\theta(M)}$ if $(y_1,\dots,y_M)\in W^{M\circ}$. \end{lemma} 
{\it Proof.} The element $\langle y_1,\dots,y_M\rangle$ is sent to 
$$\sum_{I\subseteq\{1,\dots,M\}}(-1)^{\# I}(\sum_{i\in I}y_i)^{\otimes M}
=\sum_{1\le i_1,\dots,i_M\le M}A_{i_1,\dots,i_M}y_{i_1}
\otimes\dots\otimes y_{i_M}.$$ If $S=\{1,\dots,M\}\backslash
\{i_1,\dots,i_M\}$ then $A_{i_1,\dots,i_M}=
\sum_{J\subseteq S}(-1)^{M-\# J}$, so $A_{i_1,\dots,i_M}=0$ 
if $S$ is non-empty, and $A_{i_1,\dots,i_M}=(-1)^M$ if 
$\{1,\dots,M\}=\{i_1,\dots,i_M\}$. \qed 

\begin{lemma} \label{v-samom-obxem} If $M=|\lambda|$, 
$\mu\in k$, $y_0,y_1,y_0+y_1\in W^{\circ}$ and all coordinates 
of $t_2,\dots,t_M\in W$ are algebraically independent over 
$k(y_0,y_1)$ then $$\langle y_0+y_1,t_2,\dots,t_M\rangle\otimes 
b\equiv\langle y_0,t_2,\dots,t_M\rangle\otimes b+\langle 
y_1,t_2,\dots,t_M\rangle\otimes b\bmod\ker\varphi,$$ and 
$\langle\mu y_1,t_2,\dots,t_M\rangle\otimes b\equiv\mu\langle 
y_1,t_2,\dots,t_M\rangle\otimes b\bmod\ker\varphi$. \end{lemma}
{\it Proof.} It follows from Lemmas \ref{first-rel} and \ref{calc}, 
that $\langle z_0+z_1,z_2,\dots,z_M\rangle\otimes b-
\langle z_0,z_2,\dots,z_M\rangle\otimes b-\langle z_1,\dots,z_M\rangle
\otimes b$ and $\langle\mu z_1,z_2,\dots,z_M\rangle\otimes b
-\mu\cdot\langle z_1,\dots,z_M\rangle\otimes b$ are sent to zero by 
$\varphi$, where the coordinates of $z_j$ are $x_{jm+1},\dots,x_{jm+m}$. 
As the $G$-orbits of these elements are also sent to zero by $\varphi$, 
for some $u,v\in W^{\circ}$ with coordinates algebraically independent 
over the subfield in $F$ generated over $k$ by $y_0,y_1,t_2,\dots,t_M$, 
one has the following congruences modulo the kernel of $\varphi$: 
\begin{gather}\label{per} \langle y_0+y_1,t_2,\dots,t_M
\rangle\otimes b\equiv\langle y_0+y_1+u,t_2,\dots,t_M\rangle\otimes 
b-\langle u,t_2,\dots,t_M\rangle\otimes b,\\ 
\label{vtor} \langle y_0,t_2,\dots,t_M\rangle\otimes b\equiv\langle 
y_0+u-v,t_2,\dots,t_M\rangle\otimes b-\langle u-v,t_2,\dots,t_M
\rangle\otimes b,\\ 
\label{tret}\langle y_1,t_2,\dots,t_M\rangle\otimes b\equiv\langle 
y_1+v,t_2,\dots,t_M\rangle\otimes b-\langle v,t_2,\dots,t_M\rangle
\otimes b.\end{gather} 
As $\langle y_0+y_1+u,t_2,\dots,t_M\rangle\otimes b\equiv\langle 
y_0+u-v,t_2,\dots,t_M\rangle\otimes b+\langle y_1+v,t_2,\dots,t_M\rangle
\otimes b$, and $\langle u,t_2,\dots,t_M\rangle\otimes b\equiv\langle 
u-v,t_2,\dots,t_M\rangle\otimes b+\langle v,t_2,\dots,t_M\rangle\otimes b$, 
the left hand side of the congruence (\ref{per}) is congruent to the sum 
of the left hand sides of the congruences (\ref{vtor}) and 
(\ref{tret}) modulo $\ker\varphi$. \qed 

\begin{lemma} \label{raz-obxenie} Let $(y_1,\dots,y_M)\in\{0\}\times 
W^{(M-1)\circ}\cup W^{M\circ}$ and let the coordinates of $t_{ij}\in 
W^{\circ}$ be algebraically independent over $k(y_1,\dots,y_M)$, where 
$1\le i\le M$ and $2\le j\le M$. Set $[0]:=0$ and $\langle 0,y_2,\dots,
y_M\rangle:=0$. Then \begin{multline} \label{linejnost} \langle 
y_1,\dots,y_M\rangle\otimes b\equiv\sum_{J\subseteq\{2,\dots,M\}}
(-1)^{\# J}\langle y_1,\sum_{s\in\{1\}\cup J}t_{s2},\dots,
\sum_{s\in\{1\}\cup J}t_{sM}\rangle\otimes b\\ 
-\sum_{\emptyset\neq I\subseteq\{2,\dots,M\}}(-1)^{\# I}
\langle y_1,y_2+\sum_{i\in I}t_{2i},\dots,y_M+\sum_{i\in I}t_{Mi}
\rangle\otimes b\bmod\ker\varphi.\end{multline} \end{lemma}
{\it Proof.} It follows from the identities \begin{multline*}
[\sum_{s\in J}y_s]=\sum_{\emptyset\neq I\subseteq\{2,\dots,M\}}
(-1)^{\# I}\left([\sum_{s\in J}\sum_{i\in I}t_{si}]
-[\sum_{s\in J}(y_s+\sum_{i\in I}t_{si})]\right)\\ -\langle
\sum_{s\in J}y_s,\sum_{s\in J}t_{s2},\dots,\sum_{s\in J}
t_{sM}\rangle\end{multline*} that 
\begin{multline}\langle y_1,\dots,y_M\rangle=
\sum_{\emptyset\neq I\subseteq\{2,\dots,M\}}(-1)^{\# I}
\left(\langle\sum_{i\in I}t_{1i},\dots,\sum_{i\in I}t_{Mi}\rangle
\right. \\ \left. -\langle y_1+\sum_{i\in I}t_{1i},\dots,y_M+
\sum_{i\in I}t_{Mi}\rangle\right)-\sum_{J\subseteq\{1,\dots,M\}}
(-1)^{\# J}\langle\sum_{s\in J}y_s,\sum_{s\in J}t_{s2},\dots,
\sum_{s\in J}t_{sM}\rangle.\end{multline} Then Lemma \ref{v-samom-obxem}, 
applied to the summands containing $y_1$, implies that 
\begin{multline*}\langle y_1,\dots,y_M\rangle\otimes b\equiv
\sum_{J\subseteq\{2,\dots,M\}}(-1)^{\# J}\langle y_1,
\sum_{s\in\{1\}\cup J}t_{s2},\dots,\sum_{s\in\{1\}\cup J}t_{sM}
\rangle\otimes b\\ -\sum_{\emptyset\neq I\subseteq\{2,\dots,M\}}
(-1)^{\# I}\langle y_1,y_2+\sum_{i\in I}t_{2i},\dots,y_M+
\sum_{i\in I}t_{Mi}\rangle\otimes b+\langle 0,y_2,\dots,y_M
\rangle\otimes b,\end{multline*} so we get (\ref{linejnost}). \qed

\begin{lemma} \label{voobxe} If $M=|\lambda|$, $\mu\in k$ and
$(z_j,y_2,\dots,y_M),(\sum\limits_{i=1}^Nz_i,y_2,\dots,y_M),(\mu z_1,y_2,
\dots,y_M)\in W^{M\circ}$ for all $1\le j\le N$ then \begin{equation}
\label{mnogolinej} \langle\sum_{j=1}^Nz_j,y_2,\dots,y_M\rangle\otimes 
b\equiv\sum_{j=1}^N\langle z_j,y_2,\dots,y_M\rangle\otimes b\bmod\ker
\varphi,\end{equation} and $\langle\mu z_1,y_2,\dots,y_M\rangle\otimes 
b\equiv\mu\langle z_1,y_2,\dots,y_M\rangle\otimes b\bmod\ker\varphi$.
\end{lemma}
{\it Proof.} If $N=2$ then (\ref{mnogolinej}) follows from Lemmas 
\ref{v-samom-obxem} and \ref{raz-obxenie}. If $N\ge 3$ then 
$$\langle\sum_{j=1}^Nz_j,y_2,\dots,y_M\rangle\otimes b\equiv\langle
\sum_{j=3}^Nz_j-u,y_2,\dots,y_M\rangle\otimes b+\langle z_1+z_2+u,
y_2,\dots,y_M\rangle\otimes b$$ for any sufficiently general $u\in 
W^{\circ}$. By the induction assumption, this is congruent to 
$\sum_{j=3}^N\langle z_j,y_2,\dots,y_M\rangle\otimes b-\langle u,
y_2,\dots,y_M\rangle\otimes b+\langle z_1,y_2,\dots,y_M\rangle
\otimes b+\langle z_2+u,y_2,\dots,y_M\rangle\otimes b\equiv
\sum_{j=1}^N\langle z_j,y_2,\dots,y_M\rangle\otimes b$. \qed 

\begin{lemma} \label{to-tenz} The $k$-linear map $k[W^{M\circ}]
\longrightarrow\bigotimes_k^MW$, given by $[(y_1,\dots,y_M)]\mapsto 
y_1\otimes\dots\otimes y_M$, is surjective and its kernel is spanned 
over $k$ by $[(y_0,\dots,y_{j-1}+y_j,\dots,y_M)]-[(y_0,\dots,
\widehat{y_{j-1}},\dots,y_M)]-[(y_0,\dots,\widehat{y_j},\dots,
y_M)]$ and $\mu[(y_1,\dots,y_M)]-[(y_1,\dots,\mu y_j,\dots,y_M)]$ 
for all $y_0,\dots,y_M\in W^{\circ}$ and all $\mu\in k^{\times}$. \end{lemma}
{\it Proof.} By Zorn's lemma, there exists a maximal subset $S$ in 
$W^{\circ}$ consisting of $k$-linear independent elements. If $S$ 
does not generate $W$ then the $k$-linear envelope of $S$ does not 
contain $W^{\circ}$, i.e., an element $y\in W^{\circ}$ $k$-linear 
independent over $S$, so $S\cup\{y\}$ is a bigger subset in 
$W^{\circ}$ consisting of $k$-linear independent elements. 
This contradiction shows that $S$ is a $k$-basis of $W$. 

For any $y\in W^{\circ}$ and any $z\in W$ there exist at most $m$ 
values of $\mu\in k$ such that $y+\mu z\not\in W^{\circ}$, since this 
condition is equivalent to vanishing of the $\Omega^m_{F/k}$-valued 
polynomial $(dy_1+\mu dz_1)\wedge\dots\wedge(dy_m+\mu dz_m)$ of degree 
$\le m$ in $\mu$ with non-zero constant term. Let us show that the map 
$(k^{\times}S)^M\cap W^{M\circ}\longrightarrow S^M$ given by the 
projectivization is surjective. Indeed, let $(s_1,\dots,s_M)\in S^M$. 
For all but $\le m$ values of $\mu\in k^{\times}$ one has $s_1+\mu s_2
\in W^{\circ}$. Fix one of such $\mu$ and set $s_2':=\mu s_2$. For all 
but $\le 3m$ values of $\mu\in k^{\times}$ one has $s_1+\mu s_3,s_2'+\mu 
s_3,s_1+s_2'+\mu s_3\in W^{\circ}$. Fix one of such $\mu$ and set $s_3'
:=\mu s_3$. Proceeding further this way, we get an element 
$(s_1,s'_2,\dots,s'_M)\in((k^{\times}S)^M)\cap W^{M\circ}$ projecting 
onto $(s_1,\dots,s_M)$. 

Fix a section of the projection $(k^{\times}S)^M\cap W^{M\circ}
\longrightarrow S^M$. Denote by $\widetilde{S^M}$ the image of $S^M$ 
under this section. Then $\widetilde{S^M}$ considered as a subset in 
$k[W^{M\circ}]$ maps to a basis of $\bigotimes_k^MW$, which shows 
the surjectivity. 

For the injectivity is suffices to check that $k[\widetilde{S^M}]$ 
maps onto $k[W^{M\circ}]$ modulo the relations. For an element 
$w=(w_1,\dots,w_M)\in W^{M\circ}$ set $l(w):=\sum_{j=1}^Ml_j(w)\ge M$, 
where $l_j(w)$ is the number of non-zero coordinates of $w_j$ in 
the basis $S$. 

By induction on $l(w)$ we are going to show that 
$[w]$ is in the image of $k[\widetilde{S^M}]$. 

If $l(w)=M$ then $w_j=\mu_js_j$ for all $1\le j\le M$, where 
$(s_1,\dots,s_M)\in\widetilde{S^M}$. For any sufficiently general 
$\nu_2,\dots,\nu_M\in k^{\times}$ one has \begin{multline*}
[w]\equiv\nu_2^{-1}\cdots\nu_M^{-1}[(w_1,\nu_2w_2,\dots,\nu_Mw_M)]
\equiv\mu_1\nu_2^{-1}\cdots\nu_M^{-1}[(s_1,\nu_2w_2,\dots,\nu_Mw_M)] \\ 
\equiv\mu_1\mu_2\nu_3^{-1}\cdots\nu_M^{-1}[(s_1,s_2,\nu_3w_3,\dots,\nu_M
w_M)]\equiv\dots\equiv\mu_1\cdots\mu_M[(s_1,\dots,s_M)].\end{multline*} 

The induction step: if, for instance, $l_1(w)\ge 2$ then for all but 
$\le l_1(w)m+m$ values of $\mu\in k^{\times}$ one has 
$[w]\equiv\mu^{-1}[(\mu w_1,w_2,\dots,w_M)]\equiv\mu^{-1}\sum_{s\in S}
[(\mu\mu_ss,w_2,\dots,w_M)]$, where $w_1=\sum_{s\in S}\mu_ss$ 
is a finite sum. By the induction assumption, the summands 
$[(\mu\mu_ss,w_2,\dots,w_M)]$ are in the image of $k[\widetilde{S^M}]$, 
and thus, $[w]$ is also there. \qed 

\vspace{5mm}

Let ${\mathfrak m}$ be the kernel of the multiplication
map $F\otimes_kF\stackrel{\times}{\longrightarrow}F$. 
The map $F\otimes_k(F/k)\longrightarrow{\mathfrak m}$, given 
by $\sum_jz_j\otimes\overline{y_j}\mapsto\sum_jz_j\otimes y_j
-(\sum_jz_jy_j)\otimes 1$ is clearly an isomorphism, so we can 
use the notation ${\mathfrak m}$ instead of $F\otimes_k(F/k)$, 
and the multiplicative structure of the ideal ${\mathfrak m}$.
\begin{lemma} \label{tenz-max-id-gen} The element $\alpha_q:=
(x_1\otimes 1-1\otimes x_1)^{s_1}\otimes\dots\otimes
(x_q\otimes 1-1\otimes x_q)^{s_q}\in\bigotimes^q_F{\mathfrak m}$ 
generates the sub-object ${\mathfrak m}^{s_1}
\otimes_F\dots\otimes_F{\mathfrak m}^{s_q}$. \end{lemma}
{\it Proof.} We need to show that for any collection of $\beta_i
\in{\mathfrak m}^{s_i}$ the element $\beta_1\otimes\dots\otimes 
\beta_q$ belongs to the $F[G_{F/k}]$-submodule generated by $\alpha_q$. 
Set $\alpha:=x_1\otimes 1-1\otimes x_1$. Then the $G_{F/k}$-orbit 
of $\alpha^s$ contains $(\sum_{j=1}^sa_j(y_j\otimes 1-1\otimes y_j))^s$ 
for any $a_j\in k$ and $y_j\in F$ such that $\sum_{j=1}^sa_jy_j\not\in k$. 
The $k$-span of such elements with fixed $y_1,\dots,y_s$ contains 
$\prod_{j=1}^s(y_j\otimes 1-1\otimes y_j)$. Such products generate 
${\mathfrak m}^s$ as an ideal. Moreover, they generate ${\mathfrak m}^s$ 
as a $F\otimes_kk$-module: $(1\otimes b)\prod_{j=1}^s(y_j\otimes 1-
1\otimes y_j)=((by_1\otimes 1-1\otimes by_1)-(y_1\otimes 1)
(b\otimes 1-1\otimes b))\prod_{j=2}^s(y_j\otimes 1-1\otimes y_j)=
(by_1\otimes 1-1\otimes by_1)\prod_{j=2}^s(y_j\otimes 1-1\otimes y_j)
-(y_1\otimes 1)(b\otimes 1-1\otimes b)\prod_{j=2}^s
(y_j\otimes 1-1\otimes y_j)$. 

This implies that $\beta_i=\sum_{j=1}^{s_i}f_{ij}\cdot\sigma_{ij}
\alpha^{s_i}$ for some $\sigma_{ij}\in G_{F/k}$ and $f_{ij}\in F$. 
The $G_{F/k}$-orbit of $\alpha_q$ contains $\alpha':=
(z_1\otimes 1-1\otimes z_1)^{s_1}\otimes\dots\otimes
(z_q\otimes 1-1\otimes z_q)^{s_q}$, where $z_1,\dots,z_q\in F$ are 
algebraically independent over the subfield in $F$ generated over 
$k$ by all $f_{ij},\sigma_{ij}x_1$. 

For each pair $(i,j)$ such that $1\le i\le q$ and $1\le j\le s_i$ 
there exists an element $\xi_{ij}\in G_{F/k}$ fixing all $f_{\lambda\mu},
\sigma_{\lambda\mu}x_1$ and the elements $z_{i+1},\dots,z_q$, such that 
$\xi_{ij}z_{\mu}=z_{\mu}+\sigma_{ij}x_1$. Then $(\sum_{j=1}^{s_q}f_{qj}
(\xi_{qj}-1)^{s_q})\dots(\sum_{j=1}^{s_1}f_{1j}(\xi_{1j}-1)^{s_1})
(\alpha')=\beta_1\otimes\dots\otimes\beta_q$. \qed 

\begin{corollary} \label{tenz-max-id} Any homomorphism $F\otimes_k
\bigotimes^M_k(F/k)\longrightarrow V$ factors through $\bigotimes
^M_F({\mathfrak m}/{\mathfrak m}^s)$ for some $s\ge 1$. \end{corollary}
{\it Proof.} For any integer $s\ge 1$ the element $\alpha_s:=
(x_1\otimes 1-1\otimes x_1)^s\otimes x_2\otimes\dots\otimes x_M=
\sum_{j=0}^s(-1)^j\binom{s}{j}x_1^{s-j}\otimes x_1^j\otimes x_2
\otimes\dots\otimes x_M\in({\mathfrak m}^s\otimes_k\bigotimes
^{M-1}_k(F/k))_M^{({\rm Aff}_M)_u}$ is homogeneous of degree
$s+M-1$. As $V_M$ is finite-dimensional, the image of $\alpha_s$
in $V_M$ is zero for all sufficiently big $s$. Note that $\alpha_s$
generates ${\mathfrak m}^s\otimes_k\bigotimes^{M-1}_k(F/k)$
as an $F$-semi-linear representation of $G$.

This implies that the image of $U$ is a quotient of ${\mathfrak m}/
{\mathfrak m}^s\otimes_k\bigotimes^{M-1}_k(F/k)$ for some 
$s\ge 1$, and therefore, any homomorphism 
$F\otimes_k\bigotimes^M_k(F/k)=\bigotimes^M_F{\mathfrak m}
\longrightarrow V$ factors through $\bigotimes^M_F
({\mathfrak m}/{\mathfrak m}^s)$ for some $s\ge 1$. \qed 

\begin{theorem} \label{irr-form} Any (finitely generated) object $V$ 
of ${\mathcal A}$ is a quotient of a (finite) direct sum of objects 
of type $\bigotimes^q_F({\mathfrak m}/{\mathfrak m}^s)$ for some 
$q,s\ge 1$ and $F$, if $k=\overline{{\mathbb Q}}$. In particular, 
any irreducible object of ${\mathcal A}$ is a direct summand of the 
tensor algebra $\bigotimes^{\bullet}_F\Omega^1_{F/k}$. \end{theorem}
{\it Proof.} $V$ is generated by $V_m$ for some $m\ge 0$. By Lemma 
\ref{semisimplicity}, $V_m^{({\rm Aff}_m)_u}$ is a semi-simple 
${\rm GL}_mk$-module generating $V$. As it is explained at the beginning 
of this section, $V$ is a quotient of a direct sum of $U$'s corresponding 
to irreducible direct summands of $V_m^{({\rm Aff}_m)_u}$. By 
Lemmas \ref{voobxe} and \ref{to-tenz}, $V$ is a quotient of a 
direct sum of $F\otimes_k\bigotimes^M_k(F/k)$ for some $M\ge 0$. 
Then the conclusion follows from Corollary \ref{tenz-max-id} 
and the identities ${\mathfrak m}^j/{\mathfrak m}^{j+1}=
{\rm Sym}^j_F({\mathfrak m}/{\mathfrak m}^2)$ and 
${\mathfrak m}/{\mathfrak m}^2=\Omega^1_{F/k}$. \qed

\begin{corollary} Any finitely generated object of ${\mathcal A}$ is 
of finite length, if $k=\overline{{\mathbb Q}}$. \qed \end{corollary}

\subsection{Ext's in ${\mathcal A}$} \label{proekt-pro-obra}
\begin{lemma}\label{morf-na-tenz-step} ${\rm Hom}_{{\mathcal C}}
(\bigotimes^r_F{\mathfrak m},\bigotimes^q_F({\mathfrak m}/{\mathfrak m}
^{N+1}))={\rm Hom}_{{\mathcal A}}(\bigotimes^r_F({\mathfrak m}/
{\mathfrak m}^{N+1}),\bigotimes^q_F({\mathfrak m}/{\mathfrak m}^{N+1}))$ 
admits a natural $k$-basis identified with the set $P=P(q,r,N)$ of the 
surjections $\{1,\dots,r\}\longrightarrow\{1,\dots,q\}$ with fibres of 
cardinality $\le N$, if $N\ge 1$, $q,r\ge 0$ and $q+r\ge 1$. In particular 
(take $r\ge q=1$), any subobject of ${\mathfrak m}/{\mathfrak m}^{N+1}$ 
is of type ${\mathfrak m}^r/{\mathfrak m}^{N+1}$. \end{lemma}
{\sc Example.} $P=\emptyset$ if $q>r$, or if $r>qN$; 
$\# P=q!$ if $q=r$; $\# P=1$ if $N\ge r\ge q=1$. 

{\it Proof.} By Lemma \ref{tenz-max-id-gen}, ${\mathfrak m}^{s_1}
\otimes_F\dots\otimes_F{\mathfrak m}^{s_r}$ is generated 
by the element $\otimes^r_{j=1}(x_j\otimes 1-1\otimes x_j)^{s_j}\in 
(\bigotimes^r_{j=1}{\mathfrak m}^{s_j})_r^{({\rm Aff}_r)_u}$ 
of weight $(s_1,\dots,s_r)$ with respect to $(k^{\times})^r\subseteq
{\rm GL}_rk:={\rm Aff}_r/({\rm Aff}_r)_u$. The central weights of 
$(\bigotimes^q_F({\mathfrak m}/{\mathfrak m}^{N+1}))_r^{({\rm Aff}_r)_u}$ 
are contained in the interval $[q,qN]$, so $\bigotimes^r_{j=1}
{\mathfrak m}^{s_j}$ is mapped to 0 if $\sum_{j=1}^rs_j\not\in[q,qN]$. 
In particular, the morphisms factor through $\bigotimes^r_F
({\mathfrak m}/{\mathfrak m}^{qN-r+2})$, and are zero if $r>qN$. 

The elements $\pi_{\varphi}:=\otimes^q_{i=1}\prod_{\varphi(u)=i}
(x_u\otimes 1-1\otimes x_u)$ for all surjections $\varphi\in P$ span 
the $(\underbrace{1,\dots,1}\limits_r)$-eigenspace of $(\bigotimes
^q_F({\mathfrak m}/{\mathfrak m}^{N+1}))_r^{({\rm Aff}_r)_u}$. 
Any morphism of $\bigotimes^r_F{\mathfrak m}$ is determined 
by the image $\sum_{\varphi\in P}\lambda_{\varphi}\pi_{\varphi}$ 
of the generator $\otimes^r_{j=1}(x_j\otimes 1-1\otimes x_j)$ 
for some collection of $\lambda_{\varphi}\in k$. \qed 

\begin{lemma}\label{A-spl} ${\mathcal A}$ splits as ${\mathcal V}ec_k
\oplus{\mathcal A}^{\circ}$, where ${\mathcal V}ec_k$ is the category 
of finite-dimensional $k$-vector spaces and ${\mathcal A}^{\circ}$ is
the full subcategory of ${\mathcal A}$ with objects $V$ 
such that $V^{G_{F/k}}=0$. \end{lemma}
{\it Proof.} For any $V\in{\mathcal A}$ set $V^{\circ}:=
\bigcap_{\varphi\in{\rm Hom}_{{\mathcal C}}(V,F)}\ker\varphi$. 
It follows from Theorem \ref{irr-form} and Lemma \ref{morf-na-tenz-step} 
that $V=(V^{G_{F/k}}\otimes_kF)\oplus V^{\circ}$, and ${\rm Ext}_{{\mathcal A}}
^{\ast}(F,{\mathcal A}^{\circ})={\rm Ext}_{{\mathcal A}}^{\ast}
({\mathcal A}^{\circ},F)=0$. The equivalence is given by 
$V\mapsto(V^{G_{F/k}},V^{\circ})$. \qed 

\vspace{4mm}

\label{ves-filtr}
Define the following decreasing ``weight'' filtration on the objects
$V$ of ${\mathcal A}$: $W^qV$ is the sum of the images of all morphisms
to $V$ from $\bigotimes^{\ge q}_F{\mathfrak m}$. Clearly, $W^{\bullet}$
is functorial and multiplicative. By Theorem \ref{irr-form}, $gr^q_WV$
is a finite direct sum of direct summands of $\bigotimes^q_F\Omega^1_F$.

\begin{corollary}\label{edinstv-proekt} ${\mathcal A}^{\circ}$
has no non-zero projective objects. \end{corollary}
{\it Proof.} Let $P\in{\mathcal A}^{\circ}$ be a projective object and
$\xi_2:P\longrightarrow S^{\lambda}_F\Omega^1_F$ be its irreducible 
quotient for a Young diagram $\lambda$, where $|\lambda|$ is minimal such 
that $W^{|\lambda|+1}P\neq P$. Then, for any $s\ge 2$, there is a lifting
$\xi_s:P\longrightarrow S^{\lambda}_F({\mathfrak m}/{\mathfrak m}^s)$ of
$\xi_2$. By Theorem \ref{irr-form}, there exist $q,a\ge 1$ and a morphism     
$\bigotimes_F^q({\mathfrak m}/{\mathfrak m}^a)\longrightarrow P$ such 
that its composition with $\xi_2$ is non-zero. Then its composition 
with any $\xi_s$ is also non-zero. By Lemma \ref{morf-na-tenz-step},
${\rm Hom}_{{\mathcal A}}(\bigotimes^q_F({\mathfrak m}/
{\mathfrak m}^a),S^{\lambda}_F({\mathfrak m}/{\mathfrak m}^N))=0$
for any $N\ge a+q$, leading to contradiction. \qed

\begin{lemma} \label{pre-pro} One has ${\rm Ext}^1_{{\mathcal A}}
(\bigotimes^q_F({\mathfrak m}/{\mathfrak m}^N),V)=0$ for any 
$V\in{\mathcal A}$ of finite length, $q\ge 1$ and 
$N>$ the maximal weight of $V$. \end{lemma}
{\it Proof.} Induction on the length of $V$ reduces the problem to 
the case of irreducible $V$. Let $0\longrightarrow V\longrightarrow E
\stackrel{\pi}{\longrightarrow}\bigotimes^q_F({\mathfrak m}/{\mathfrak m}^N)
\longrightarrow 0$ be an extension. By Theorem \ref{irr-form}, there 
is a surjection of a direct sum of objects of type $\bigotimes^p_F
({\mathfrak m}/{\mathfrak m}^a)$ onto $E$. By Lemma 
\ref{morf-na-tenz-step}, ${\rm Hom}_{{\mathcal A}}(\bigotimes^{\neq q}
_F({\mathfrak m}/{\mathfrak m}^a),\bigotimes^q_F({\mathfrak m}/
{\mathfrak m}^2))=0$, so there is a morphism of a direct sum of objects 
of type $\bigotimes^q_F({\mathfrak m}/{\mathfrak m}^a)$ to $E$ 
surjective over $\bigotimes^q_F({\mathfrak m}/{\mathfrak m}^2)$. 
As the latter is semi-simple, there is a morphism of $\bigoplus
_{|\lambda|=q}S^{\lambda}_F({\mathfrak m}/{\mathfrak m}^a)$ to $E$ 
surjective over $\bigotimes^q_F({\mathfrak m}/{\mathfrak m}^2)$. 
By Lemma \ref{morf-na-tenz-step}, its composition with $\pi$ is 
surjective, and therefore, the weights of its kernel are $\ge N$, 
so it does not intersect $V$. In other words, the extension splits. \qed 

\begin{corollary}\label{4.12} The following pro-representable functor
on ${\mathcal A}$ $${\rm Hom}_{{\mathcal C}}({\mathfrak m}^{s_1}
\otimes_F\dots\otimes_F{\mathfrak m}^{s_q},-)=\lim_{\longrightarrow}
{\rm Hom}_{{\mathcal A}}(({\mathfrak m}^{s_1}/{\mathfrak m}^N)\otimes_F
\dots\otimes_F({\mathfrak m}^{s_q}/{\mathfrak m}^N),-)$$ is exact if 
and only if $s_1=\dots=s_q=1$. \end{corollary}
{\it Proof.} Let $V\longrightarrow V'$ be a surjection in ${\mathcal A}$ 
and $\xi:\bigotimes^q_F{\mathfrak m}\longrightarrow V'$ be a morphism in 
${\mathcal C}$. We have to show that $\xi$  factors through $V$. By Lemma 
\ref{tenz-max-id-gen}, the image of $\xi$ is cyclic. Let $V''$ be the 
cyclic sub-object of $V$ generated by a pre-image of a generator of the 
image of $\xi$. Then the kernel $K$ of $V''\longrightarrow{\rm Im}(\xi)$ 
is of finite length. As $\xi$ factors through $\bigotimes^q_F
({\mathfrak m}/{\mathfrak m}^N)$ for some $N\gg 0$, and Lemma \ref{pre-pro} 
implies that ${\rm Ext}^1_{{\mathcal A}}(\bigotimes^q_F({\mathfrak m}/
{\mathfrak m}^N),K)=0$, $\xi$ factors through $V$. 

The rest follows from the fact that the projection ${\mathfrak m}^s
\longrightarrow{\mathfrak m}^s/{\mathfrak m}^{N+s}$ does not lift 
to ${\mathfrak m}^s\longrightarrow\bigotimes^s_F({\mathfrak m}/
{\mathfrak m}^{N+1})$, if $s\ge 2$: neither non-zero morphism 
$\bigotimes^s_F{\mathfrak m}\longrightarrow\bigotimes^s_F({\mathfrak m}
/{\mathfrak m}^{N+1})$ factors through ${\mathfrak m}^s$, if $N\ge 2$. \qed

\begin{corollary}\label{fin-dim-ext} If $V\in{\mathcal A}$ is of 
finite type then $\dim_k{\rm Ext}_{{\mathcal A}}^j(V,V')<\infty$ 
for any $j\ge 0$ and any $V'\in{\mathcal A}$. If $V\in{\mathcal A}$ 
is irreducible and ${\rm Ext}_{{\mathcal A}}^1({\mathfrak m}/
{\mathfrak m}^q,V)\neq 0$ for some $q\ge 2$ then $V\cong
{\mathfrak m}^q/{\mathfrak m}^{q+1}$ and ${\rm Ext}_{{\mathcal A}}
^1({\mathfrak m}/{\mathfrak m}^q,V)\cong k$. \end{corollary}
{\it Proof.} If $V\in{\mathcal A}$ is of finite type then, by Theorem 
\ref{irr-form}, it admits a resolution $\dots\longrightarrow P_2
\longrightarrow P_1\longrightarrow P_0$ whose terms are finite direct
sums of objects of type $\bigotimes^s_F{\mathfrak m}$. By Lemma 
\ref{tenz-max-id-gen}, the terms of the complex ${\rm Hom}_{{\mathcal C}}
(P_{\bullet},V')$ are finite-dimensional over $k$ and, by Corollary 
\ref{4.12}, it calculates ${\rm Ext}_{{\mathcal A}}^{\bullet}(V,V')$. \qed 

\begin{corollary} The filtration $W^{\bullet}$ is 
strictly compatible with the surjections. \end{corollary}
{\it Proof.} Let $V\longrightarrow V'$ be a surjection in ${\mathcal A}$. 
Then, by Corollary \ref{4.12}, any morphism $\bigotimes^q_F
{\mathfrak m}\longrightarrow V'$ factors through $V$. \qed 

\section{``Coherent'' sheaves in smooth topology}
\label{smooth-topol}
Let ${\mathfrak S}m_k$ be the category of locally dominant morphisms 
of smooth $k$-schemes. Consider on ${\mathfrak S}m_k$ the (pre-)topology, 
where the covers are surjective smooth morphisms. Clearly, the covers 
are stable under the base changes. 

By definition, the structure presheaf ${\mathcal O}$ of 
${\mathfrak S}m_k$ associates to any $Y\in{\mathfrak S}m_k$ 
its $k$-algebra of regular functions ${\mathcal O}(Y)$. 
Clearly, ${\mathcal O}$ is a sheaf in this topology.

A sheaf ${\mathcal F}$ on ${\mathfrak S}m_k$ is 
``{\sl (quasi-)coherent}'' if its values ${\mathcal F}(Y)$ are 
endowed with ${\mathcal O}(Y)$-module structures and its restriction 
to the small \'etale site of $Y$ is a (quasi-)coherent sheaf for any 
$Y\in{\mathfrak S}m_k$. 

\begin{lemma}\label{4.14} Let $X\longrightarrow Y$ be an \'etale morphism 
of smooth varieties over $k$ sending a point $q\in X$ to a point $p\in Y$. 
Then ${\mathfrak m}^s_q/{\mathfrak m}^N_q={\mathcal O}_q\otimes
_{{\mathcal O}_p}({\mathfrak m}^s_p/{\mathfrak m}^N_p)$ for any $s\le N$, 
where ${\mathfrak m}_q:=\ker({\mathcal O}_q\otimes_k{\mathcal O}_q
\stackrel{\times}{\longrightarrow}{\mathcal O}_q)$. \end{lemma}
{\it Proof.} One has ${\mathfrak m}_q/{\mathfrak m}^2_q={\mathcal O}_q
\otimes_{{\mathcal O}_p}({\mathfrak m}_p/{\mathfrak m}^2_p)$, so applying 
${\rm Sym}^s_{{\mathcal O}_q}$ we get ${\mathfrak m}^s_q/{\mathfrak m}
^{s+1}_q={\rm Sym}^s_{{\mathcal O}_q}({\mathfrak m}_q/{\mathfrak m}^2_q)
={\mathcal O}_q\otimes_{{\mathcal O}_p}{\rm Sym}^s_{{\mathcal O}_p}
({\mathfrak m}_p/{\mathfrak m}^2_p)={\mathcal O}_q\otimes_{{\mathcal O}_p}
({\mathfrak m}^s_p/{\mathfrak m}^{s+1}_p)$. 
The induction on $N-s$ gives the conclusion:
$$\begin{array}{ccccccccc} 0 & \to & {\mathfrak m}^{s+1}_q/
{\mathfrak m}^N_q & \to & {\mathfrak m}^s_q/{\mathfrak m}^N_q 
& \to & {\mathfrak m}^s_q/{\mathfrak m}^{s+1}_q & \to & 0 \\
 && \| && \bigcup && \| &&  \\
0 & \to & {\mathcal O}_q\otimes_{{\mathcal O}_p}({\mathfrak m}^{s+1}_p
/{\mathfrak m}^N_p) & \to & {\mathcal O}_q\otimes_{{\mathcal O}_p}
({\mathfrak m}^s_p/{\mathfrak m}^N_p) & \to &
{\mathcal O}_q\otimes_{{\mathcal O}_p}({\mathfrak m}^s_p
/{\mathfrak m}^{s+1}_p) & \to & 0 \end{array}$$
\qed

\begin{corollary} \label{funk-A-presh} The category ${\mathcal A}$ 
is equivalent to the category of ``coherent'' sheaves on 
${\mathfrak S}m_k$, if $k=\overline{{\mathbb Q}}$. \end{corollary}
{\it Proof.} Fix an embedding over $k$ of the function field of each 
connected component of each smooth $k$-variety into $F$. Then, for any 
$V\in{\mathcal A}$, $Y\in{\mathfrak S}m_k$ and a point $q\in Y$ define 
an ${\mathcal O}_q$-lattice ${\mathcal V}_q\subset V^{G_{F/k(Y)}}$ as 
follows. Let ${\mathcal O}_p\subseteq{\mathcal O}_q$ be an \'etale 
extension of a local subring in $F$ of a closed point $p$ of a 
projective space. 

Any object $V$ of ${\mathcal A}$ is a quotient of a direct sum
of objects of type $\bigotimes^s_F({\mathfrak m}/{\mathfrak m}^N)$.
Then, as it is true for 
$\bigotimes^s_F({\mathfrak m}/{\mathfrak m}^N)$ (Lemma \ref{4.14}), 
it follows that the module ${\mathcal V}_p\subset V$ provided by the 
exact functor ${\mathcal S}$, cf. \S\ref{no-1}, is independent of the 
choice of the projective space, and ${\mathcal V}_q:={\mathcal O}_q
\otimes_{{\mathcal O}_p}{\mathcal V}_p\subset V$ is independent of 
${\mathcal O}_p$.

This determines a locally free coherent sheaf ${\mathcal V}_Y$ on $Y$ 
with the generic fibre $V^{G_{F/k(Y)}}$. 

It follows also that, for any dominant morphism $X\stackrel
{\pi}{\longrightarrow}Y$ of smooth $k$-varieties, the inclusion 
of the generic fibres $k(X)\otimes_{k(Y)}V^{G_{F/k(Y)}}\subseteq 
V^{G_{F/k(X)}}$ induces an injection of the coherent sheaves 
$\pi^{\ast}{\mathcal V}_Y\hookrightarrow{\mathcal V}_X$ on $X$, 
which is an isomorphism if $\pi$ is {\'e}tale. 

To check that ${\mathcal V}$ is a sheaf on ${\mathfrak S}m_k$, we 
need to show that for any surjective smooth morphism $X\longrightarrow 
Y$ the sequence $0\longrightarrow{\mathcal V}(Y)
\stackrel{\beta}{\longrightarrow}{\mathcal V}(X)
\stackrel{p^{\ast}_1-p^{\ast}_2}{\longrightarrow}
{\mathcal V}(X\times_YX)$ is exact. As ${\mathcal V}_X$ is a sheaf 
in Zariski topology on $X$, it suffices to treat the case of affine 
$X$ and $Y$. In the case $V=\bigotimes^s_F({\mathfrak m}/
{\mathfrak m}^N)$, which is sufficient by Theorem \ref{irr-form} 
and Lemma \ref{morf-na-tenz-step}, this amounts to the exactness 
of the sequence $0\longrightarrow\bigotimes^s_B({\mathfrak m}_B/
{\mathfrak m}_B^N)\longrightarrow\bigotimes^s_A({\mathfrak m}_A/
{\mathfrak m}_A^N)\longrightarrow\bigotimes^s_{A\otimes_BA}
({\mathfrak m}_{A,B}/{\mathfrak m}_{A,B}^N)$, where 
$B$ is a smooth $k$-algebra of finite type, $A$ is a smooth 
$B$-algebra of finite type, ${\mathfrak m}_C:=\ker(C\otimes_kC
\stackrel{\times}{\longrightarrow}C)$ for any $k$-algebra $C$, and 
${\mathfrak m}_{A,B}:={\mathfrak m}_{A\otimes_BA}$. But this is clear. 

Conversely, a ``coherent'' sheaf ${\mathcal V}$ on ${\mathfrak S}m_k$ 
is sent to the object $\lim\limits_{\longrightarrow}{\mathcal V}(U)$, 
where $U$ runs over the spectra of regular subalgebras in $F$ of 
finite type over $k$. (As $F$ is the union of its regular subalgebras 
of finite type over $k$, $\lim\limits_{\longrightarrow}{\mathcal V}(U)$ 
is an $(F=\lim\limits_{\longrightarrow}{\mathcal O}(U))$-module. The 
action of an element $\sigma\in G$ comes as the limit of isomorphisms 
$\sigma^{\ast}:{\mathcal V}(U)\stackrel{\sim}{\longrightarrow}
{\mathcal V}(U')$, where $U={\bf Spec}(A)$ and $U={\bf Spec}(\sigma(A))$ 
induced by the isomorphism $U'\stackrel{\sim}{\longrightarrow}U$.) \qed 

\begin{lemma}\label{birac-invar} For any ``quasi-coherent'' flat (as 
${\mathcal O}$-module) sheaf ${\mathcal V}$ on ${\mathfrak S}m_k$ 
the $k$-space ${\mathcal V}(Y)$ is a birational invariant of proper 
$Y$. If ${\mathcal V}$ is ``coherent'' then ${\mathcal V}({Y'})$ 
generates the (generic fibre of the) sheaf ${\mathcal V}_{Y'}$ 
for appropriate finite covers $Y'$ of $Y$. \end{lemma}
{\it Proof.} According to Hironaka, for any pair of smooth proper 
birational $k$-varieties $Y,Y''$ there is a smooth proper $k$-variety 
$Y'$ and birational $k$-morphisms $Y'\stackrel{\pi}{\longrightarrow}Y$ 
and $Y'\longrightarrow Y''$. Let $Z\subset Y$ be the subset consisting 
of points $z$ such that $\pi:\pi^{-1}(z)\to z$ is not an isomorphism. 
It is a subvariety of codimension $\ge 2$. As ${\mathcal V}$ is 
torsion-free, one has ${\mathcal V}(Y)\longrightarrow{\mathcal V}(Y')
\stackrel{i^{\ast}}{\hookrightarrow}{\mathcal V}(U)$,\footnote{To show 
that $i^{\ast}$ is also injective, choose an affine covering $\{U_j\}$ 
of $Y'$, and a dense affine subset $U'\subseteq U$. As sum of ample 
divisors is ample, any intersection of open affine subsets is again 
affine, so $\{U_j\cap U'\}$ is an affine covering of $U'$. Then the 
diagram $$\begin{array}{ccccc} && {\mathcal V}(Y')&\hookrightarrow 
& \bigoplus_i{\mathcal V}(U_i)\\ & \lefteqn{i^{\ast}}\swarrow & 
\downarrow && \phantom{\varphi}\downarrow\varphi \\ {\mathcal V}(U) 
& \rightarrow & {\mathcal V}(U')&\hookrightarrow & \bigoplus_i
{\mathcal V}(U_i\cap U')\end{array}$$ is commutative, and $\varphi$ 
is injective since ${\mathcal V}(U_i\cap U')={\mathcal O}(U_i\cap U')
\otimes_{{\mathcal O}(U_i)}{\mathcal V}(U_i)$ and ${\mathcal V}$ 
is torsion-free.} where $U:=Y-Z\stackrel{i}{\hookrightarrow}Y'$ is the 
section of $\pi$. It suffices to check that for any affine $Y$ one has 
${\mathcal V}(Y)\stackrel{\sim}{\longrightarrow}{\mathcal V}(U)$. 
Choose an affine covering $\{U_j\}$ of $U$. Then $0\longrightarrow
{\mathcal V}(U)\longrightarrow\bigoplus_j{\mathcal O}(U_j)
\otimes_{{\mathcal O}(Y)}{\mathcal V}(Y)\longrightarrow\bigoplus_{i,j}
{\mathcal O}(U_i\cap U_j)\otimes_{{\mathcal O}(Y)}{\mathcal V}(Y)$ 
is exact, so, as $0\longrightarrow{\mathcal O}(U)={\mathcal O}(Y)
\longrightarrow\bigoplus_j{\mathcal O}(U_j)\longrightarrow
\bigoplus_{i,j}{\mathcal O}(U_i\cap U_j)$ is also exact, 
we get ${\mathcal V}(Y)={\mathcal V}(U)$. \qed 

\vspace{4mm}

{\sc Remark.} If ${\mathcal V}:Y\mapsto\Omega^j_{k(Y)}/\Omega^j(Y)$ then 
the sequence $0\longrightarrow\Omega^j(Y)\longrightarrow\Omega^j_{k(Y)}
\longrightarrow{\mathcal V}(Y)\longrightarrow H^1(Y,\Omega^j_Y)
\longrightarrow 0$ is exact, so ${\mathcal V}(Y)$ is birationally 
invariant if and only if $j=0$: for any closed smooth $Z\subset
{\mathbb P}^{j+1}=Y$ of codimension 2 such that $\Omega^{j-1}(Z)\neq 0$ 
one has $H^1(Y',\Omega^j_{Y'})\cong H^1(Y,\Omega^j_Y)\oplus
\Omega^{j-1}(Z)$, where $Y'$ is the blow-up of $Y$ along $Z$. 

\vspace{4mm}

Then, using Lemma \ref{birac-invar}, we get a left exact (non faithful) 
functor (with faithful restriction to the subcategory of ``coherent'' 
sheaves) $$\left\{\mbox{flat ``quasi-coherent'' sheaves 
on ${\mathfrak S}m_k$}\right\}\stackrel{\Gamma}{\longrightarrow}\left
\{\mbox{smooth representations of $G_{F/k}$ over $k$}\right\}$$ 
given by ${\mathcal V}\mapsto\lim\limits_{\longrightarrow}
\Gamma(Y,{\mathcal V}_Y)$, where $Y$ runs over the smooth proper 
models of subfields in $F$ of finite type over $k$. This functor 
is not full, and the objects in its image are highly reducible, 
e.g., $\Gamma(\Omega^1_{/k})\cong\bigoplus_A(A(F)/A(k))\otimes
_{{\rm End}(A)}\Gamma(A,\Omega^1_{A/k})$, where $A$ runs over the 
set of isogeny classes of simple abelian varieties over $k$. If 
${\mathcal V}$ is ``coherent'' and $\Gamma(Y,{\mathcal V}_Y)$ has the 
Galois descent property then $\Gamma({\mathcal V})$ is admissible. 
However, there is no Galois descent property in general.

{\sc Example.} Let $Y'$ be a smooth projective hyperelliptic curve
$y^2=P(x)$, considered as a 2-fold cover of the projective line $Y$.
Then, for ${\mathcal V}_Y=(\Omega^1_{Y/k})^{\otimes 2}$, the
section $y^{-2}(dx)^2=P(x)^{-1}(dx)^2$ is a Galois invariant
element of $\Gamma(Y',{\mathcal V}_{Y'})$,
which is not in $\Gamma(Y,{\mathcal V}_Y)=0$.

\section{${\mathcal A}/{\mathcal A}_{>m}$}
\label{faktor-kat}
The only finite-dimensional objects of
${\mathcal A}$ are direct sums of copies of $F$, so the
category ${\mathcal A}$ is far from being tannakian.
However, ${\mathcal A}$ admits a decreasing filtration
by Serre subcategories ${\mathcal A}_{>m}$ such that all
${\mathcal A}/{\mathcal A}_{>m}$ are again abelian tensor
categories and their objects are finite-dimensional.
The category ${\mathcal A}/{\mathcal A}_{>m}$ is not rigid.

Let ${\mathcal A}_{>m}$ be the full subcategory of ${\mathcal A}$
with objects $V$ such that $V_m=0$. Clearly, ${\mathcal A}_{>m}$ is
a Serre subcategory of ${\mathcal A}$. Moreover, it is an ``ideal''
in ${\mathcal A}$ in the sense that the tensor product functor
${\mathcal A}_{>m}\times{\mathcal A}\longrightarrow{\mathcal A}$
factors through ${\mathcal A}_{>m}$, so the quotient abelian
category ${\mathcal A}/{\mathcal A}_{>m}$ carries a tensor structure.

By definition, the objects of ${\mathcal A}/{\mathcal A}_{>m}$ are
the objects of ${\mathcal A}$, but the morphisms are defined by
${\rm Hom}_{{\mathcal A}/{\mathcal A}_{>m}}(V,V')={\rm Hom}
_{{\mathcal A}}(\langle V_m\rangle,V'/(V')_{>m})=
{\rm Hom}_{{\mathcal A}/{\mathcal A}_{>m}}(V,\langle V'_m\rangle)$,
where $\langle V_m\rangle$ denotes the semi-linear subrepresentation
of $V$ generated by $V_m$ and $(V')_{>m}$ is the maximal subobject
of $V'$ in ${\mathcal A}_{>m}$.\footnote{The functor ${\mathcal A}
\longrightarrow{\mathcal A}_{>m}$, $V\mapsto(V)_{>m}$ is right adjoint
to inclusion functor ${\mathcal A}_{>m}\longrightarrow{\mathcal A}$.
In particular, it is left exact.} In particular, $V\cong\langle
V_m\rangle$ in ${\mathcal A}/{\mathcal A}_{>m}$.

{\sc Example.} ${\mathcal A}/{\mathcal A}_{>0}$ is equivalent
to the category of finite-dimensional $k$-vector spaces.

The functor ${\mathcal A}/{\mathcal A}_{>m}\longrightarrow
\mathfrak{SL}^u_m$, $V\mapsto V_m$ is exact, faithful and tensor.
Note also that the objects of ${\mathcal A}/{\mathcal A}_{>m}$ are
{\sl finite-dimensional}. Namely, $\bigwedge^{\dim_{K_m}V_m+1}V=0$.

Let $\Phi$ be a monoid of one-dimensional objects of ${\mathcal A}
/{\mathcal A}_{>m}$, such as $(\Omega^m_{F/k})^{\otimes N}$
for any $N\ge 0$. The set $\Phi$ is partially ordered:
$\omega\le\eta$ if there is $\xi\in{\mathcal A}/{\mathcal A}_{>m}$
such that $\eta\cong\omega\otimes\xi$. In particular,
$\omega\le\omega\otimes\eta$ and $\eta\le\omega\otimes\eta$.
If $k=\overline{{\mathbb Q}}$ then $\Phi$ consists of some
(symmetric) $F$-tensor powers of $\Omega^m_{F/k}$.
\begin{lemma} The $k$-vector space ${\rm Hom}_{{\mathcal A}/
{\mathcal A}_{>m}}(V\otimes\omega,V'\otimes\omega)$ is
finite-dimensional and independent of $\omega\in\Phi$ 
for $\omega$ sufficiently big. \end{lemma}
{\it Proof.} For any $\omega,\eta\in\Phi$ such that $\omega\le\eta$
(i.e., $\eta\cong\omega\otimes\xi$ for some $\xi\in{\mathcal A}/
{\mathcal A}_{>m}$) the twist by $\xi$ defines a canonical inclusion
${\rm Hom}_{{\mathcal A}/{\mathcal A}_{>m}}(V\otimes\omega,
V'\otimes\omega)\subseteq{\rm Hom}_{{\mathcal A}/{\mathcal A}_{>m}}
(V\otimes\eta,V'\otimes\eta)$. The $k$-vector spaces
\begin{multline*}{\rm Hom}_{{\mathcal A}/{\mathcal A}_{>m}}(V,V')=
{\rm Hom}_{{\mathcal A}}(\langle V_m\rangle,V'/(V')_{>m})
\longrightarrow{\rm Hom}_{{\mathcal A}/{\mathcal A}_{>m}}(V\otimes U,
V'\otimes U)\\ ={\rm Hom}_{{\mathcal A}}(\langle V_m\otimes_{K_m}U_m
\rangle,(V'\otimes U)/(V'\otimes U)_{>m})\subseteq{\rm Hom}_{K_m
\langle G_{K_m/k}\rangle}(V_m\otimes_{K_m}U_m,V'_m\otimes_{K_m}U_m)
\end{multline*} are finite-dimensional. On the other hand,
\begin{multline*}{\rm Hom}_{{\mathcal A}}(\langle V_m\rangle,
V'/(V')_{>m})\subseteq{\rm Hom}_{K_m\langle G_{K_m/k}\rangle}
(V_m,V'_m)\\ \subseteq{\rm Hom}_{K_m\langle G_{K_m/k}\rangle}
(V_m\otimes_{K_m}U_m,V'_m\otimes_{K_m}U_m)\end{multline*}
for any $U\in{\mathcal C}$ with $U_m\neq 0$,
where the second equality takes place if and only if
$\dim_{K_m}U_m=1$, e.g., for $U\in\Phi$. \qed

\vspace{5mm}

Let $Ob({\mathcal A}^+_{\Phi,m}):=Ob({\mathcal A})$ and ${\rm Hom}
_{{\mathcal A}^+_{\Phi,m}}(V,V'):={\rm Hom}_{{\mathcal A}/
{\mathcal A}_{>m}}(V\otimes\omega,V'\otimes\omega)$ for
sufficiently big $\omega\in\Phi$. Then $\otimes\omega:
{\mathcal A}^+_{\Phi,m}\longrightarrow{\mathcal A}^+_{\Phi,m}$
is a fully faithful functor, so we can invert objects in $\Phi$
to get a category ${\mathcal A}_{\Phi,m}:={\mathcal A}^+_{\Phi,m}
[\Phi^{-1}]$. If $\Phi$ is the set of all one-dimensional objects
of ${\mathcal A}/{\mathcal A}_{>m}$ then
${\mathcal A}_m:={\mathcal A}_{\Phi,m}$ is tannakian.

\vspace{5mm}

\noindent
{\sl Acknowledgement.} {\small I would like to thank Uwe Jannsen for 
many inspiring discussions, and in particular for suggesting (in 
Spring 2003) that smooth representations can be considered as 
sheaves in some topology. 
I am indebted to Dmitry Kaledin for thorough reading of a previous 
version of this paper and proposing numerous clarification. 
I am grateful to the Max-Planck-Institut 
f\"ur Mathematik in Bonn and to Regensburg University for 
their hospitality, and to the Max-Planck-Institut and to 
Alexander von Humboldt-Stiftung for the financial support. 

}

\vspace{5mm}

$$\begin{array}{l} \mbox{Independent University of Moscow} \\
\mbox{121002 Moscow} \\ \mbox{B.Vlasievsky Per. 11} 
\\ \mbox{{\tt marat@mccme.ru}} \end{array}\quad\mbox{and}\quad
\begin{array}{l}\mbox{Institute for Information} \\ 
\mbox{Transmission Problems} \\ 
\mbox{of Russian Academy of Sciences} \end{array}$$
\end{document}